\input amstex
\documentstyle{amsppt}
\document
\topmatter
\title
The rigid analytical regulator and Drinfeld modular curves
\endtitle
\title
The rigid analytical regulator and $K_2$ of Drinfeld modular curves
\endtitle
\author
Ambrus P\'al
\endauthor
\date
May 31, 2009.
\enddate
\address Department of Mathematics, 180 Queen's Gate, Imperial College,
London SW7 2AZ, United Kingdom\endaddress
\email a.pal\@imperial.ac.uk\endemail
\abstract We  evaluate a rigid analytical analogue of the 
Beilinson-Bloch-Deligne regulator on certain explicit elements in the $K_2$ of Drinfeld modular curves, constructed from analogues of modular units, and relate its value to special values of $L$-series using the Rankin-Selberg method.
\endabstract
\footnote" "{\it 2000 Mathematics Subject Classification. \rm 11G05
(primary), 11G40 (secondary).}
\endtopmatter

\heading 1. Introduction
\endheading

\definition{Motivation 1.1} In the paper [10] which is now classical B. Gross formulated a generalization of his original $p$-adic analogue of Stark's conjecture in a form which makes good sense both over number fields and function fields. This conjecture was proved by D. Hayes for function fields in [12]. In this paper Hayes gave an explicit rigid analytical construction of Stark units and expressed them in terms of special values of $L$-functions
using this explicit construction. This paper is part of the project to formulate and prove results which generalize Hayes's theorem the same way as Beilinson's conjectures generalize Stark's. In a previous paper ([22]) we constructed a rigid analytical regulator analogous to the classical Beilinson-Bloch-Deligne regulator refining the tame regulator in case of Mumford curves. In our current work we express the value of this regulator on certain explicit elements of the $K_2$ group of Drinfeld modular curves, which are analogues of A.
Beilinson's construction using modular units, in terms of special values of $L$-functions.
Using the function field analogue of the Shimura-Taniyama-Weil conjecture we derive a formula for every elliptic curve defined over the rational
function field of transcendence degree one a finite field having split
multiplicative reduction at the point at infinity analogous to the classical theorem of Beilinson on the $K_2$ of elliptic curves defined over the rational
number field. 
\enddefinition
In the rest of this introductory chapter we first describe the rigid analytic regulator for Tate elliptic curves then define the $\infty$-adic $L$-function of elliptic curves of the type mentioned above and formulate our main theorem.
\definition{Notation 1.2} Let $F_{\infty}$ be a field complete with respect to a discrete valuation and let $\Cal O_{\infty}$ be its valuation ring. There is a canonical way to extend the absolute value of $F_{\infty}$ induced by its valuation to its algebraic closure. Let $\Bbb C_{\infty}$ denote the completion of the algebraic closure of $F_{\infty}$ with respect to this absolute value and let $|\cdot|$ denote the absolute value induced by the completion process.
Let $|\Bbb C_{\infty}|$ denote the set of values of the latter. Let $\Bbb P^1$ denote the projective line over $\Bbb C_{\infty}$. We call a set $D\subset\Bbb P^1$ an open disc if it is the image of the set $\{z\in\Bbb C_{\infty}||z|<1\}$ under a M\"obius transformation. Recall that a subset $U$ of $\Bbb P^1$ is a connected rational subdomain, if it is non-empty and it is the complement of the union of finitely many pair-wise disjoint open discs. Let $\partial U$ denote the set of these complementary open discs. Let $\Cal O(U)$ and $\Cal O^*(U)$ denote the algebra of holomorphic functions on $U$ and the group of invertible elements of this algebra, respectively.  For each $f\in\Cal O(U)$ let $\|f\|$ denote $\sup_{z\in U}|f(z)|$. This is a finite number, and makes $\Cal O(U)$ a Banach algebra over $\Bbb C_{\infty}$. The latter is the closure of the subalgebra of restrictions of
rational functions with respect to the supremum norm $\|\cdot\|$ by definition. For every real number $0<\epsilon<1$ we define the sets
$\Cal O_{\epsilon}(U)=\{f\in\Cal O(U)|\|1-f\|\leq\epsilon\}$, and
$\Cal U_{\epsilon}=\{z\in\Bbb C_{\infty}||1-z|\leq\epsilon\}$. Recall that a function $f:\Bbb C^*_{\infty}\rightarrow\Bbb C_{\infty}$ is holomorphic if its restriction $f|_U$ is holomorphic for every connected rational subdomain $U\subset\Bbb C^*_{\infty}$. For every $x\in\Bbb P^1$ and every pair of rational non-zero functions $f$, $g\in\Bbb C_{\infty}((t))$ on the projective line let $\{f,g\}_x$ denote the tame symbol of the pair $(f,g)$ at $x$. Let $\Cal M(\Bbb C^*_{\infty})$ denote the field of meromorphic functions of $\Bbb C^*_{\infty}$. For every field $L$ let $K_2(L)$ denote the Milnor $K_2$ of the field
$L$. Finally for every $x\in\Bbb C_{\infty}$ and positive number $\rho\in|\Bbb C_{\infty}|$ let $D(x,\rho)$ denote the open disc $\{z\in\Bbb C_{\infty}||z-x|<\rho\}$. The following result is an immediate consequence of the results of [22]. 
\enddefinition
\proclaim{Theorem 1.3} For every $0<r\in|\Bbb C_{\infty}|$ there is a unique homomorphism:
$$\{\cdot\}_r:K_2(\Cal M(\Bbb C^*_{\infty}))\rightarrow\Bbb C^*_{\infty}$$
with the following properties:
\roster
\item"$(i)$" for every pair of  rational functions $f$, $g\in\Cal M(\Bbb C^*_{\infty})^*$ we have:
$$\{f\otimes g\}_r=\prod_{x\in D(0,r)}\{f,g\}_x,$$
\item"$(ii)$" for every real number $0<\epsilon<1$ and functions $f\in\Cal M(\Bbb C^*_{\infty})\cap\Cal O_{\epsilon}(U)$ and $g\in\Cal M(\Bbb C^*_{\infty})\cap\Cal O^*(U)$ we have $\{f,g\}_r\in \Cal U_{\epsilon}$ where $U$ is a connected rational subdomain $U\subset\Bbb C^*_{\infty}$ such that $D(0,r)\in\partial U$.\ $\square$
\endroster
\endproclaim
\definition{Notation 1.4} For every field $K$, for any variety $V$ defined over $K$ and for any extension $L$ of $K$ let $V_L$ denote the base change of $V$ to $L$. For every field $K$ and regular irreducible projective curve $C$ defined over $K$ let $\Cal F(C)$ denote the function field of the curve $C$ over $K$. For every closed point $x$ of $C$ there is a tame symbol at $x$ which is a homomorphism from $K_2(\Cal F(C))$ into the multiplicative group of the residue field at $x$. We define the group $K_2(C)$ as the intersection of the kernels of all tame symbols. (In this paper
we will sometimes use the somewhat incorrect notation $K_2(X)$ to denote $H^2_{\Cal M}(X,\Bbb Z(2))$ for various types of spaces $X$ as the latter is rather awkward.) Let $E$ be an elliptic curve defined over $F_{\infty}$ which has a rigid-analytic Tate uniformization over $F_{\infty}$. The latter is equivalent to the property that the special fiber of the N\'eron model of $E$ over the spectrum of $\Cal O_{\infty}$ is split multiplicative.  Let $\theta:\Bbb C^*_{\infty}\rightarrow E(\Bbb C_{\infty})$ be the Tate uniformization (over $\Bbb C_{\infty}$). The latter induces a homomorphism
$$\theta^*:\Cal F(E_{F_{\infty}})\rightarrow\Cal M(\Bbb C^*_{\infty})$$
by pull-back which in turn induces a homomorphism $K_2(\Cal F(E_{F_{\infty}}))\rightarrow K_2(\Cal M(\Bbb C^*_{\infty}))$ which will be denoted by the same symbol by slight abuse of notation. 
\enddefinition
\proclaim{Proposition 1.5} For every $k\in K_2(E_{F_{\infty}})$ and $0<r\in|\Bbb C_{\infty}|$ we have $\{\theta^*(k)\}_r\in F_{\infty}^*$ and the latter is independent of the choice of $r$.\ $\square$
\endproclaim
Let $\{\cdot\}:K_2(E_{F_{\infty}})\rightarrow F^*_{\infty}$ denote the homomorphism defined by the common value of the regulators $\{\theta^*(\cdot)\}_r$. 
\definition{Definition 1.6} For every field $K$ let $\overline K$ denote its separable closure. Let $F$ denote the function field of $X$, where the latter is a geometrically connected smooth projective curve defined over the finite field $\Bbb F_q$ of characteristic $p$. Fix a closed point $\infty$ of the curve $X$ and let $E$ be an elliptic curve defined over $F$ which has split multiplicative reduction at $\infty$. For every closed point $x$ of $X$ let $\deg(x)$ and $L_x(E,t)$ denote the degree of $x$ and the local factor of the Hasse-Weil $L$-function of $E$ at $x$, respectively. The latter is an element of $\in\Bbb Z[[t]]$. Let $\psi^*_E(x^n)\in\Bbb Z$ denote the unique number such that
$$L_x(E,t)=\sum_{n=0}^{\infty}\psi^*_E(x^n)t^{n\deg(x)}.$$
Let $K$ be a number field and let $\Delta$ denote its ring of integers. Let $\chi:\text{Gal}(\overline F|F)\rightarrow K^*$ be a $K$-valued one-dimensional Galois representation of $F$ which has finite image. Note that $\chi$ is automatically almost everywhere unramified and its image lies in $\Delta^*$. Let $\Gamma$ denote the quotient of $\text{Gal}(\overline F|F)$ by the kernel of $\chi$. Assume that $\chi$ splits at $\infty$ and  let $\goth m$ be an effective divisor whose support does not contain $\infty$ and the conductor of $\chi$ and $E$ divides $\goth m$ and $\goth m\infty$, respectively. (Note that such an $\goth m$ exists because we assumed that $E$ has split multiplicative reduction at $\infty$.) For every Galois group $G$ of a finite abelian extension of $F$ and for every closed point $x$ of $X$ where $G$ is unramified let $\phi^G_x$ denote the image of a geometric Frobenius at $x$ in $G$. The element $\phi^G_x\in G$ is well-defined as $G$ is abelian. Assume now that $G$ is the Galois group of a finite abelian extension of $F$ which only ramifies at $\infty$. We define the $L$-function $\Cal L^G_{\goth m}(E,\chi,t)$ as the Euler product:
$$\Cal L^G_{\goth m}(E,\chi,t)=
\prod_{x\notin supp(\goth m\infty)}\left(\sum_{n=0}^{\infty}\psi^*_E(x^n)
\chi(\phi^{\Gamma}_x)(\phi^G_x)^nt^{n\deg(x)}\right)\in\Delta[G][[t]],$$
where $supp(\goth d)$ denotes the support of any effective divisor
$\goth d$ on $X$. The infinite product $\Cal L^G_{\goth m}(E,\chi,t)$ is well-defined, as the constant term of every factor appearing in the product is $1$, and there are only finitely many factors with a term of degree less than $m$ for any positive integer $m$. Actually even more is true:
\enddefinition
\proclaim{Proposition 1.7} The power series $\Cal L^G_{\goth m}(E,\chi,t)$ is an element of $\Delta[G][t]$.
\endproclaim
\definition{Definition 1.8} An important consequence of the proposition above is that the polynomial $\Cal L^G_{\goth m}(E,\chi,t)$ can be evaluated at $1$, i.e. the element $\Cal L^G_{\goth m}(E,\chi,1)\in\Bbb\Delta[G]$ is well-defined. Let $G_{\infty}$ denote the Galois group of the maximal abelian extension of $F$ unramified at every closed point $x$ of $X$ different from $\infty$. It is a profinite group. Also note that if $H$ denotes the Galois group of the maximal abelian extension of $F$ unramified at every closed point $x$ of $X$ and totally split at $\infty$, then the kernel of the natural projection $G_{\infty}\rightarrow H$ is canonically isomorphic to the profinite completion of $F_{\infty}^*/\Bbb F^*_q$, the multiplicative group of the completion $F_{\infty}$ of $F$ with with respect to the valuation at $\infty$ divided out by the multiplicative group of the constant field of $X$. (Note that this notation is compatible with what we have introduced in 1.2 and 1.4.) For any ring $R$ and abelian profinite group $M$ let 
 $R[[M]]$ denote the $R$-dual of the ring of continuous functions $f:M\rightarrow R$, where $f$ is continuous with respect to the discrete topology on $R$ and the Krull topology on $M$. The ring $R[[M]]$ is also the projective limit of $R$-coefficient group rings of the finite quotients of $M$. The elements $\Cal L^G_{\goth m}(E,\chi,1)$ satisfy the obvious compatibility, so their limit defines an element $\Cal L_{\goth m}(E,\chi)$ in $\Delta[[G_{\infty}]]$, which we will call the $\infty$-adic $L$-function of $E$ twisted with $\chi$. For every $M$ as above let $I_M\triangleleft\Delta[[M]]$ denote the kernel of the natural augmentation map $\Delta[[M]]\rightarrow\Delta$. We will usually drop the subscript $M$ to ease notation. It is known that the group $I_M/I_M^2$ is naturally isomorphic to $M\otimes\Delta$. Finally let $\theta'\in M\otimes\Delta$ denote the class of any $\theta\in I_M$ in $I_M/I_M^2$.
\enddefinition
\proclaim{Proposition 1.9} We have $\Cal L_{\goth m}(E,\chi)\in I$ and
$\Cal L_{\goth m}(E,\chi)'\in F_{\infty}^*/\Bbb F^*_q\otimes\Delta$.
\endproclaim
Let $L$ denote the Galois extension of $F$ whose Galois group is
$\Gamma$. By our assumptions the field $L$ has an imbedding into $F_{\infty}$ which extends the canonical inclusion $F\subset F_{\infty}$. Fix once and for all such an imbedding. By slight abuse of notation let $\{\cdot\}:K_2(E_L)\otimes K\rightarrow F^*_{\infty}\otimes K$ denote also the composition of the homomorphism $K_2(E_L)\otimes K\rightarrow K_2(E_{F_{\infty}})\otimes K$ induced by the imbedding above and the unique $K$-linear extension of the homomorphism $\{\cdot\}$. Assume that $F=\Bbb F_q(T)$ is the rational function field of transcendence degree one over $\Bbb F_q$, where $T$ is an indeterminate, and $\infty$ is the point at infinity on $X=\Bbb P^1_{\Bbb F_q}$. Also assume that $\chi$ is non-trivial. Now we are able to state our main result:
\proclaim{Theorem 1.10} There is an element $\kappa_E(\chi)\in
K_2(E_L)\otimes K$ such that
$$\{\kappa_E(\chi)\}=L(E,q^{-1})\Cal L_{\goth m}(E,\chi)'
\text{\ in\ }F_{\infty}^*\otimes K.$$
\endproclaim
It is easy to deduce that the valuation of $\Cal L_{\goth m}(E,\chi)'$ with respect to $\infty\otimes\text{id}_K$ is equal to $-L_{\goth m}(E,\chi,1)$ from the interpolation property. (For the explanation of this notation see the next chapter.) Deligne's purity theorem implies that the latter is non-zero under mild, purely local conditions on $\chi$ and $\goth m$. If the special value $L(E,q^{-1})$ also happens to be non-zero we get that the element $\kappa_E(\chi)\in K_2(E_L)\otimes K$ is not torsion hence our main result is non-vacuous. 
\definition{Contents 1.11} In the next chapter we prove the basic properties
of the $L$-function $\Cal L_{\goth m}(E,\chi)$ by simple cohomological means. We introduce our mail tool, which we call double Eisenstein series, in the
third chapter. They are really analogous to the product of two Eisenstein
series in the classical setting, but they cannot be written as such due to
the lack of logarithm in positive characteristic. Here we also establish their basic properties, among them Proposition 3.5, which is analogous to analytic continuation. The link between double Eisenstein series and the rigid analytic regulator of elements in $K_2$ analogous to Beilinson's construction
is provided by the Krokecker limit formula 4.10 of the fourth chapter. The fifth chapter is somewhat technical: it identifies function field analogues of modular units with the rigid analytic functions appearing in the previous chapter and studies the action of the Hecke algebra on the
source and target groups of the rigid analytic regulator.
We execute the principal calculation of the paper in the sixth
chapter. Perhaps the crucial reason why the Rankin-Selberg computation can
be carried out is that the double Eisenstein series do become a product of two series after the first step of the calculation. In the seventh chapter we use the function field analogue of the Shimura-Taniyama-Weil conjecture as well as its explicit description due to Gekeler and Reversat to conclude the proof of our main result. The aim of the last chapter is to prove a useful lemma on the action of correspondences on motivic cohomology groups which is used in the fifth chapter. 
\enddefinition
\definition{Acknowledgment 1.12} I wish to thank the CRM and the IH\'ES for
their warm hospitality and the pleasant environment they created for
productive research, where this article was written.
\enddefinition

\heading 2. The $\infty$-adic $L$-functions of elliptic curves
\endheading

\definition{Definition 2.1} Note that for a finite group $G$ we have
$\Delta[[G]]=\Delta[G]$ naturally. Let $M$ be an abelian profinite group, let $H$ be a finite quotient of $M$ and let $K$ denote the kernel of the quotient map $M\rightarrow H$. We let $I^M_H$ denote the ideal of the quotient map $\Delta[[M]]\rightarrow\Delta[H]$. It is obvious that the augmentation ideal $I=I^M_{\{1\}}$ and $I^M_H\subseteq I$ for any $H$.
\enddefinition
\proclaim{Lemma 2.2} We have $\theta'\in K\otimes\Delta$ for any
$\theta\in I^M_H$.
\endproclaim
\definition{Proof} The same as the proof of Lemma 3.9 of [21] whose claim is just slightly
different.\ $\square$
\enddefinition
\definition{Notation 2.3} Let $E$ be an elliptic curve defined over $F$ which has split
multiplicative reduction at $\infty$ as in the introduction whose notation we are going to use
without further notice. Let $G$ be the Galois group of a finite abelian extension of $F$ which
only ramifies at $\infty$ and let $H(G)$ denote the maximal quotient of $G$ such that the
corresponding abelian extension of $F$ is unramified at every closed point $x$ of $X$ and
totally split at $\infty$.
\enddefinition
\proclaim{Proposition 2.4} The following holds:
\roster
\item"$(i)$" the power series $\Cal L^G_{\goth m}(E,\chi,t)$ is an element of $\Delta[G][t]$,
\item"$(ii)$" we have $\Cal L^G_{\goth m}(E,\chi,1)\in I^G_{H(G)}$.
\endroster
\endproclaim
\definition{Proof}  Let $l$ be a prime different from $p$. The Gal$(\overline F|F)$-module $H^1(E_{\overline F},\Bbb Q_l)$ is absolutely irreducible because the curve $E$ is not isotrivial. Let $\rho$ denote the corresponding $l$-adic Galois representation. For every character $\phi:G\rightarrow\overline{\Bbb Q}_l^*$
let the same symbol denote the corresponding homomorphism $\overline{\Bbb Q}_l[G][[t]]\rightarrow\overline{\Bbb Q}_l[[t]]$ and the corresponding $l$-adic Galois representation by the usual abuse of notation. For every $l$-adic Galois representation $\psi$ which is unramified at almost all closed points of $X$ will use the same symbol to denote the constructible $l$-adic sheaf on $X$ which is the direct image of $\psi$ with respect to the generic point Spec$(F)\rightarrow X$. Fix an imbedding of $K$ into $\overline{\Bbb Q}_l$. This way the Galois representation $\chi$ becomes an $l$-adic representation, too. The series $\Cal L^G_{\goth m}(E,\chi,t)\in\overline{\Bbb Q}_l[[t]]$ is characterized by the property:
$$\phi(\Cal L^G_{\goth m}(E,\chi,t))=L(X(\goth m\infty),
\rho\otimes\chi\phi,t)$$
for every character $\phi:G\rightarrow\overline{\Bbb Q}_l^*$ where 
where $X(\goth d)$ denotes complement of the support of any effective divisor $\goth d$ in $X$ and
$L(U,\psi,t)$ denotes the Grothendieck $L$-function of any constructible $l$-adic sheaf $\psi$ on a variety $U$ over $\Bbb F_q$. The $l$-adic Galois representation $\rho\otimes\chi\phi$ is absolutely
irreducible, therefore the twisted $L$-function $L(X(\goth m\infty),
\rho\otimes\chi\phi,t)$ is a polynomial for every character $\phi:G\rightarrow\overline{\Bbb Q}_l^*$ by the Grothendieck-Verdier
formula. Hence so does $\Cal L^G_{\goth m}(E,\chi,t)$ as claim $(i)$ says. For every character $\phi:H(G)\rightarrow\overline{\Bbb Q}_l^*$ let the same symbol denote the composition of the quotient map $G\rightarrow H(G)$ and the character $\phi$ as well. In this case the restriction of the $l$-adic Galois representation corresponding to
$\phi$ to the decomposition group at $\infty$ is trivial. The same holds for $\chi$ by assumption. Moreover $E$ has split multiplicative reduction at $\infty$ so we have:
$$\phi(\Cal L^G_{\goth m}(E,\chi,t))=(1-t^{\deg(\infty)})L(X(\goth m),
\rho\otimes\chi\phi,t)$$
for every such character. As the twisted $L$-function $L(X(\goth m),\rho\otimes\chi\phi,t)$ is a polynomial by the Grothendieck-Verdier formula, we have $\phi(\Cal L^G_{\goth m}(E,\chi,1))=0$ for every such character
as well. The latter is equivalent to the property that $\Cal L^{H(G)}_{\goth m}(E,\chi,1)$ is zero as claim
$(ii)$ says.\ $\square$
\enddefinition
As we explained in Definition 1.8 part $(i)$ of the proposition above implies that the object $\Cal L_{\goth m}(E,\chi)$ is well-defined. For every group $M$ let $\widehat M$ denote its profinite
completion and let $\infty:\widehat{F^*_{\infty}/\Bbb F^*_q}
\otimes\Delta\rightarrow\widehat{\Delta}=\widehat{\Bbb Z}
\otimes\Delta$ denote the profinite completion of the valuation $\infty$ as well. The following proposition takes care of Proposition 1.9 and the remark after Theorem 1.10. For the sake of simple notation let $L_{\goth m}(E,\chi,t)$ denote $L(X(\goth m),\rho\otimes\chi,t)$. 
\proclaim{Proposition 2.5} The following holds:
\roster
\item"$(i)$" we have $\Cal L_{\goth m}(E,\chi)\in I$ and
$\Cal L_{\goth m}(E,\chi)'\in F_{\infty}^*/\Bbb F^*_q\otimes\Delta$.
\item"$(ii)$" we have $\infty(\Cal L_{\goth m}(E,\chi)')=-L_{\goth m}(E,\chi,1)$.
\endroster
\endproclaim
\definition{Proof} The first half of claim $(i)$ and the fact that
$\Cal L_{\goth m}(E,\chi)'\in\widehat{F^*_{\infty}/\Bbb F^*_q}\otimes
\Delta$ follows at once from claim $(ii)$ of Proposition 2.4 and Lemma 2.2 by taking the limit.
On the hand note that $F_{\infty}^*/\Bbb F^*_q\otimes\Delta$ is the pre-image of $\Delta$ with respect to $\infty$ in $\widehat{F^*_{\infty}/\Bbb F^*_q}\otimes\Delta$ hence the second half of claim $(i)$ is an immediate consequence of claim $(ii)$. Now we only have to show the latter. The profinite group $G_{\infty}$ surjects onto the
Galois group of the maximal constant field extension of $F$ which is isomorphic to $\widehat{\Bbb Z}$. This induces a surjection $\Delta[[G_{\infty}]]\rightarrow\Delta[[\widehat{\Bbb Z}]]$. The choice of a topological generator of $\widehat{\Bbb Z}$, or equivalently the choice of a system of generators of the finite quotients of $\widehat{\Bbb Z}$ compatible with the projections furnishes an injection $\Delta[t]\rightarrow\Delta[[\widehat{\Bbb Z}]]$ such that the image of $t$ is the generator. In case of the natural choice of the global geometric Frobenius as a
topological generator, the image $\phi_x$ of a geometric Frobenius at $x$ in $G_{\infty}$ maps to $t^{\deg(x)}$ for every closed point $x$ on $X$ under the map above. Hence the image of $\Cal L_{\goth m}(E,\chi)$ under this map is $\widetilde L_{\goth m}(E,\chi,t)=(1-t^{\deg(\infty)})L_{\goth m}(E,\chi,t)$ as we saw in the proof of Proposition 2.4. The ideal $I\triangleleft\Delta[[G_{\infty}]]$ maps into the augmentation ideal $J\triangleleft\Delta[[\widehat{\Bbb Z}]]$ corresponding to the trivial quotient of $\widehat{\Bbb Z}$, and the induced map $I/I^2\rightarrow J/J^2$ is the tensor product of the surjection $G_{\infty}\rightarrow\widehat{\Bbb Z}$ introduced above and the identity of $\Delta$. Since the intersection $J\cap\Delta[t]$ is the ideal generated by $t-1$, the image of $\Cal L_{\goth m}(E,\chi)'$ under the map $I/I^2\rightarrow J/J^2$ is just the derivative $\widetilde L(E,1)'\in\Delta\subset\widehat{\Delta}$. On the other hand the restriction of the surjection $G_{\infty}\rightarrow\widehat{\Bbb Z}$ to $F_{\infty}^*$ is $\deg(\infty)$ times the valuation map $\infty:F^*_{\infty}\rightarrow\Bbb Z$, so:
$$\deg(\infty)\infty(\Cal L_{\goth m}(E,\chi)')=((1-t^{\deg(\infty)})L_{\goth m}(E,\chi,t))'\bigg|_{t=1}=-\deg(\infty)L_{\goth m}(E,\chi,1)$$
as we claimed.\ $\square$
\enddefinition

\heading 3. Double Eisenstein series
\endheading

\definition{Notation 3.1} Let $|X|$, $\Bbb A$, $\Cal O$ denote set of closed points of $X$,
the ring of adeles of $F$ and its maximal compact subring of $\Bbb A$, respectively. As in the
introduction we will fix a closed point $\infty$ in the set $|X|$. For every divisor $\goth m$
of $X$ let $\goth m$ also denote the $\Cal O$-module in the ring $\Bbb A$ generated by the
ideles whose divisor is $\goth m$, by abuse of notation. For every idele $m\in\Bbb A^*$ let
the same symbol also denote the divisor of $m$ if this notation does not cause confusion. For
any closed point $v$ in $|X|$ we will let $F_v$, $\bold f_v$ and $\Cal O_v$ denote the
corresponding completion of $F$, its constant field, and its discrete valuation ring,
respectively. For every $v\in|X|$ let $v:F_v^*\rightarrow\Bbb Z$ denote the valuation
normalized such that $v(\pi_v)=\deg(v)$ for every uniformizer $\pi_v\in F_v$.
For any idele, adele, adele-valued matrix or
function defined on the above which decomposes as an infinite product
of functions defined on the individual components the subscript $v$
will denote the $v$-th component. Let $\Bbb A_f$, $\Cal O_f$ denote the restricted direct
product $\prod'_{x\neq\infty} F_x$ and the direct product $\prod_{x\neq\infty}\Cal O_x$, respectively. The
former is also called the ring of finite adeles of $F$ and the latter is its maximal compact
subring. For every $g\in GL_2(\Bbb A)$ (or $g\in\Bbb A$, $etc$.) let $g_f$ denote its finite
component in $GL_2(\Bbb A_f)$. We will consider $\Bbb A^*_f$ as well as $F^*_v$ (for every
place $v\in|X|$) as a subgroup of $\Bbb A^*$ in the natural way. Similarly we will consider
$\Bbb A_f$ and $F_v$ as a subring of $\Bbb A$ and $GL_2(\Bbb A)$ and $GL_2(F_v)$ as a subgroup of $GL_2(\Bbb A)$. 
Let $|\cdot|$ denote the normalized absolute value on the ring $\Bbb A$ and
for any idele or divisor $y$ let $\deg(y)$ denote its degree related to the
normalized absolute value by the formula $|y|=q^{-\deg(y)}$.
In accordance with our convention $|\cdot|$ will denote the absolute value with respect to
$\infty$ if its argument is in
$F_{\infty}$. For each $(u,v)\in F^2_{\infty}$ let $\|(u,v)\|$,
$\infty(u,v)$ denote $\max(|u|,|v|)$ and $\min(\infty(u),\infty(v))$,
respectively. Let $Z$ denote the center of the group scheme $GL_2$, let 
$$\Gamma_{\infty}=\left\{\left(\matrix a&b\\
c&d\endmatrix\right)\in GL_2(\Cal O_{\infty})|
\infty(c)>0\right\}$$
be the Iwahori subgroup of $GL_2(F_{\infty})$ and let
$$\Bbb K(\goth m)=\{g\in GL_2(\Cal O)|g\equiv
I\text{\ mod }\goth m\},$$
for every effective divisor $\goth m$ where $I$ is the identity matrix. We will adopt the convention which assigns $0$ or $1$ as value to the empty sum or product, respectively.
\enddefinition
\definition{Definition 3.2} Let $F^2_<$ denote the set:
$F^2_<=\{(a,b)\in F_{\infty}^2||a|<|b|\}$. Let $\goth m$
be an effective divisor on $X$ whose support does not contain $\infty$.
Let the same symbol also denote the ideal $\goth m\cap\Cal O_f$ by abuse
of notation. For every $g\in GL_2(\Bbb A)$, $(\alpha,\beta)\in(\Cal
O_f/\goth m)^2$, and $n$ integer let
$$\split W_{\goth m}(\alpha,\beta,g,n)=&\{0\neq f\in F^2|fg_f\in(\alpha,\beta)+
\goth m\Cal O_f^2,-n=\infty(fg_{\infty})\},\\
V_{\goth m}(\alpha,\beta,g,n)=&\{f\in W_{\goth m}(\alpha,\beta,g,n)|
fg_{\infty}\in F^2_<\}\text{\ and}\\
U_{\goth m}(\alpha,\beta,g,n)=&W_{\goth m}(\alpha,\beta,g,n)-V_{\goth m}
(\alpha,\beta,g,n).\endsplit$$
Also let
$$W_{\goth m}(\alpha,\beta,g_f)=\bigcup_{n\in\Bbb Z}
W_{\goth m}(\alpha,\beta,g,n),$$
$$U_{\goth m}(\alpha,\beta,g)=\bigcup_{n\in\Bbb Z}U_{\goth
m}(\alpha,\beta,g,n)\text{\ and\ }V_{\goth m}(\alpha,\beta,g)=
\bigcup_{n\in\Bbb Z}V_{\goth m}(\alpha,\beta,g,n).$$
Obviously the first set is well-defined. For every finite quotient $G$ of $F^*\backslash\Bbb A^*/\Cal O_f^*$ let $\cdot^G:\Bbb A^*\rightarrow G$ denote the quotient map. Let $E^G_{\goth m}(\alpha,\beta,\gamma,\delta,g,x,y)$ denote the $\Bbb Z[G][[x,y]](x^{-1},y^{-1})$-valued function
$$\split E^G_{\goth m}(\alpha,\beta,\gamma,\delta,g&,x,y)=\\
&{\det(g_f^{-1})^G\over
(xy)^{\deg(\det(g))}}\cdot
\sum\Sb(a,b)\in U_{\goth m}(\alpha,\beta,g)\\(c,d)\in
V_{\goth m}(\gamma,\delta,g)\endSb\!\!\!\!\!\!\!\!\!\!\!\!\det
\left(\matrix a&b\\c&d\endmatrix\right)_{\infty}^G\!\!\!
x^{2\infty((a,b)g_{\infty})}y^{\infty(2(c,d)g_{\infty})},\endsplit$$
for every $g\in GL_2(\Bbb A)$, variables $x$, $y$,  and pairs $(\alpha,\beta)$ and $(\gamma,\delta)$ as above. In order to see that this function is indeed well-defined first note that
$$\left(\matrix a&b\\c&d\endmatrix\right)=\left(\matrix a_1&b_1\\c_1&d_1\endmatrix\right)
\cdot\det(g_{\infty})^{-1}=(a_1d_1-b_1c_1)\cdot\det(g_{\infty})^{-1}$$
is non-zero where $(a_1,b_1)=(a,b)g_{\infty}$ and $(c_1,d_1)=(c,d)g_{\infty}$ because $|a_1|\geq|b_1|$ and $|c_1|<|d_1|$ by the definition of the sets $U_{\goth m}(\alpha,\beta,g)$ and $V_{\goth m}(\gamma,\delta,g)$ therefore
$$|a_1d_1-b_1c_1|=|a_1d_1|\neq0.$$
Hence the terms of the infinite sum above are defined. The sum itself is well-defined 
and $\Bbb Z[G][[x,y]](x^{-1},y^{-1})$-valued as the cardinality of the sets $U_{\goth m}(\alpha,\beta,g)$ and $V_{\goth m}(\gamma,\delta,g)$ are finite for all $n$ and zero for $n$ sufficiently small.
\enddefinition
\proclaim{Proposition 3.3} The following holds:
\roster
\item"$(i)$" the function $E^G_{\goth m}(\alpha,\beta,\gamma,\delta,g,x,y)$ is
left-invariant with respect to $GL_2(F)$ and right-invariant with respect to $\Bbb K(\goth m\infty)\Gamma_{\infty}Z(F_{\infty})$,
\item"$(ii)$" the $\Bbb C[G]$-valued infinite sum
$E^G_{\goth m}(\alpha,\beta,\gamma,\delta,g,q^{-s},q^{-t})$
converges absolutely, if $\text{\rm Re}(s)>1$ and $\text{\rm Re}(t)>1$, for every $g$.
\endroster
\endproclaim
\definition{Proof} We are going to prove claim $(i)$ first. Since  for every $\rho\in GL_2(F)$ and $n\in\Bbb Z$ we have:
$$U_{\goth m}(\alpha,\beta,\rho g,n)=U_{\goth m}
(\alpha,\beta,g,n)\rho^{-1}
\text{\ and\ }V_{\goth m}(\gamma,\delta,\rho g,n)=V_{\goth m}
(\gamma,\delta,g,n)\rho^{-1},$$
we get that
$$\split E^G_{\goth m}(\alpha,\beta,\gamma,\delta,\rho g,x,y)&=
{\det(\rho_f^{-1})^G\det(g_f^{-1})^G\over(xy)^{\deg(\det(\rho))+\deg(\det(g))}}
\\\cdot&\!\!\!\!\!\!\!\!\!\!\!\!\!\!
\sum\Sb(a,b)\in U_{\goth m}(\alpha,\beta,g)\\(c,d)\in
V_{\goth m}(\gamma,\delta,g)\endSb\!\!\!\!\!\!\!\!\!\!\!\det
\left(\matrix a&b\\c&d\endmatrix\right)_{\infty}^G\!\!\!\!\!
\det(\rho_{\infty}^{-1})^Gx^{2\infty((a,b)\rho^{-1}\rho g_{\infty})}y^
{2\infty((c,d)\rho^{-1}\rho g_{\infty})}\\
&=E^G_{\goth m}(\alpha,\beta,\gamma,\delta,g,x,y),\endsplit$$
because $\det(\rho^{-1})\in F^*$ and $\deg(\det(\rho))=0$ as the degree of every principal
divisor is zero. On the other hand for every $\lambda\in GL_2(F_{\infty})$ the set
$\{f\in F^2_{\infty}|f\lambda\in F^2_<\}$ is obviously left invariant by
$\Gamma_{\infty}Z(F_{\infty})$ hence
$$U_{\goth m}(\alpha,\beta,g\rho)=U_{\goth m}(\alpha,\beta,g)
\text{\ and\ }
V_{\goth m}(\gamma,\delta,g\rho)=V_{\goth m}(\gamma,\delta,g)$$
for every $\rho\in\Bbb K(\goth m\infty)\Gamma_{\infty}Z(F_{\infty})$ and $g\in GL_2(\Bbb A)$.
Therefore
$$\split E^G_{\goth m}(\alpha,\beta,\gamma,\delta,g\rho,x,y)&=
{\det(\rho_f^{-1})^G\det(g_f^{-1})^G\over(xy)^{\deg(\det(z))+\deg(\det(g))}}
\\\cdot&\!\!\!\!\!\!\!\!\!\!\!\!\!\!
\sum\Sb(a,b)\in U_{\goth m}(\alpha,\beta,g)\\(c,d)\in
V_{\goth m}(\gamma,\delta,g)\endSb\!\!\!\!\!\!\!\!\!\!\!\!\det
\left(\matrix a&\!\!b\\c&\!\!d\endmatrix\right)_{\infty}^G\!\!\!\!\!
x^{2\infty((a,b)g_{\infty})+\infty(\det(z))}
y^{2\infty((c,d)g_{\infty})+\infty(\det(z))}\\
&=E^G_{\goth m}(\alpha,\beta,\gamma,\delta,g,x,y),\endsplit$$
where $\rho=\kappa z$ with $\kappa\in\Bbb K(\goth m\infty)\Gamma_{\infty}$ and
$z\in Z(F_{\infty})$ because $\kappa_{\infty}$ is an isometry with respect to the norm
$\|\cdot\|$, $\deg(\det(\kappa))=0$ and $\det(\kappa_f)^G=1$ by definition. Our proof of claim
$(ii)$ is the same as the argument that may be found in [20]. The coefficient of each element
of $G$ in the series $E^G_{\goth m}(\alpha,\beta,\gamma,\delta,g,
q^{-s},q^{-t})$ is majorated by the product $E(g,s)E(g,t)$ where:
$$E(g,s)=|\det(g)|^s\!\!\!\!\!
\sum\Sb f\in F^2-\{0\}\\fg\in\Cal O_f^2\endSb\!\!\!
\parallel\!(fg)_{\infty}\!\parallel^{-2s},$$
so it will be sufficient to prove that $E(g,s)$ converges absolutely for
each $g\in GL_2(\Bbb A)$ if $\text{\rm Re}(s)>1$. For every $g\in GL_2
(\Bbb A)$ let $\Cal E(g)$ denote the sheaf on $X$ whose group of
sections is for every open subset $U\subseteq X$ is
$$\Cal E(g)(U)=\{f\in F^2|fg\in\Cal O_v^2,\forall v\in|U|\},$$
where we denote the set of closed points of $U$ by $|U|$. The sheaf $\Cal
E(g)$ is a coherent locally free sheaf of rank two. If $\Cal F_n$
denote the sheaf $\Cal F\otimes\Cal O_X(\infty)^n$ for every coherent
sheaf $\Cal F$ on $X$ and integer $n$, then for every $g\in GL_2(\Bbb A)$
and $s\in\Bbb C$ the series above can be rewritten as
$$E(g,s)=\sum_{n\in\Bbb Z}|H^0(X,\Cal E(g)_n)-H^0(X,\Cal E(g)_{n-1})|
q^{-s\deg(\Cal E(g)_n)}.$$
By the Riemann-Roch theorem for curves:
$$\dim H^0(X,\Cal F)-\dim H^0(X,K_X\otimes\Cal F^{\vee})=2-2g(X)+\deg(
\Cal F)$$
for any coherent locally free sheaf of rank two $\Cal F$ on $X$,
where $K_X$, $\Cal F^{\vee}$ and $g(X)$ is the canonical bundle on $X$,
the dual of $\Cal F$, and the genus of $X$, respectively. Because $\dim
H^0(X,\Cal F_{-n})=0$ for $n$ sufficiently large depending on $\Cal F$, we
have that
$$|H^0(X,\Cal E(g)_n)|=q^{2-2g(X)+\deg(\Cal E(g))+2n\deg(\infty)}\text{\ 
and\ }|H^0(X,\Cal E(g)_{-n})|=1,$$
if $n$ is a sufficiently large positive number. Hence
$$E(g,s)=p(q^{-s})+q^{2-2g(X)+(1-s)\deg(\Cal E(g))}
(1-q^{-\deg(\infty)})\sum_{n=0}^{\infty} q^{2n(1-s)\deg(\infty)},$$
where $p$ is a polynomial. The claim now follows from the convergence of
the geometric series.\ $\square$
\enddefinition
\definition{Definition 3.4} For every abelian group $M$ and for every finite set $S$ let $M[S]$ and $M[S]_0$ denote the group of functions $f:S\rightarrow M$ and its subgroup consisting of functions $f\in M[S]$ with the property
$$\sum_{\alpha\in S}f(\alpha)=0,$$
respectively. Let $\Cal V_{\goth m}$ denote the set $(\Cal O_f/\goth m)^2-\{0,0\}$ and for every $C\in R[\Cal V_{\goth m}]$ and $D\in R[\Cal V_{\goth m}]$ let $E^G_{\goth m}(C,D,g,x,y)$ denote the function:
$$E^G_{\goth m}(C,D,g,x,y)=\sum\Sb(\alpha,\beta)\in\Cal V_{\goth m}
\\(\gamma,\delta)\in\Cal V_{\goth m}\endSb C(\alpha,\beta)D(\gamma,\delta)
E^G_{\goth m}(\alpha,\beta,\gamma,\delta,g,x,y),$$
where $R\supseteq\Bbb Z$ is an arbitrary commutative ring. 
\enddefinition
\proclaim{Proposition 3.5} For every $C$,
$D\in R[\Cal V_{\goth m}]_0$ the function $E^G_{\goth m}(C,D,g,x,y)$ takes values in $R[G][x,y,x^{-1},y^{-1}]$.
\endproclaim
\definition{Proof} We may assume by bilinearity that $C=(\alpha,\beta)-(\gamma,\delta)$  and $D=(\epsilon,\iota)-(\kappa,\lambda)$ for some pairs $(\alpha,\beta)$, $(\gamma,\delta)$, $(\epsilon,\iota)$ and $(\kappa,\lambda)\in\Cal V_{\goth m}$. 
Pick two elements $(r,s)\in U_{\goth m}(\alpha-\gamma,\beta-\delta,g)$ and $(u,v)\in V_{\goth m}(\epsilon-\kappa,\iota-\lambda,g)$. Then for every sufficiently large natural number $n$ we have:
$$U_{\goth m}(\alpha,\beta,g,n)=\{(a+r,b+s)|(a,b)\in
U_{\goth m}(\gamma,\delta,g,n)\}\text{\ and}$$
$$V_{\goth m}(\epsilon,\iota,g,n)=\{(a+u,b+v)|(a,b)\in
V_{\goth m}(\kappa,\lambda,g,n)\}.$$
Therefore
$$\split E^G_{\goth m}(C,D,g,x,y)=&P(x,y)\\
&+{\det(g_f^{-1})^G\over(xy)^{\deg(\det(g))}}\!\!\!\!\!\!\!
\sum\Sb(a,b)\in U_{\goth m}(\gamma,\delta,g_f)\\(c,d)\in
V_{\goth m}(\kappa,\lambda,g)\endSb\!
\left[\matrix a&\!\!b\\c&\!\!d\endmatrix\right]
x^{2\infty((a,b)g_{\infty})}y^{2\infty((c,d)g_{\infty})},\endsplit$$
where $P(x,y)\in\Bbb Z[G][x,y,x^{-1},y^{-1}]$ and
$$\split\left[\matrix a&b\\c&d\endmatrix\right]=&
\det\left(\matrix a+r&b+s\\c+u&d+v\endmatrix\right)_{\infty}^G-
\det\left(\matrix a+r&b+s\\c&d\endmatrix\right)_{\infty}^G\\&-
\det\left(\matrix a&b\\c+u&d+v\endmatrix\right)_{\infty}^G+
\det\left(\matrix a&b\\c&d\endmatrix\right)_{\infty}^G.\endsplit$$
In order to finish the proof it is enough to show that the determinants
in the expression above can be paired in such a way that in every pair the
determinants have different signs and they represent the same element
in $G$ if $\max(\|(a,b)g_{\infty}\|,\|(c,d)g_{\infty}\|)$ is
sufficiently large. This follows the lemma below or its pair which we get by switching the rows of the matrices depending on whether $\|(a,b)g_{\infty}\|$ or
$\|(c,d)g_{\infty}\|$ is the larger one among the two, respectively.\ $\square$
\enddefinition
\proclaim{Lemma 3.6} For each $k\in GL_2(F_{\infty})$ let $\|k\|$ denote the maximum of absolute values of the entries of $k$. Then for every $g\in GL_2(F_{\infty})$ and $(r,s)$, $(a,b)$, $(c,d)\in F^2_{\infty}$ such that $(a,b)g\notin F^2_<$ and $(c,d)g\in
 F^2_<$ we have:
$$\left|1-\det(\left(\matrix a+r&b+s\\c&d\endmatrix\right)\cdot
\left(\matrix a&b\\c&d\endmatrix\right)^{-1})\right|
\leq{\|(r,s)\|\|g^{-1}\||\det(g)|
\over\|(a,b)g\|}.$$
\endproclaim
\definition{Proof} Using Cramer's rule we get that
$$\left(\matrix a+r&b+s\\c&d\endmatrix\right)\cdot
\left(\matrix a&b\\c&d\endmatrix\right)^{-1}=
\left(\matrix1+{rd-sc\over ad-bc}&{-rb+sa\over ad-bc}
\\0&1\endmatrix\right),$$
so
$$\left|1-\det(\left(\matrix a+r&b+s\\c&d\endmatrix\right)\cdot
\left(\matrix a&b\\c&d\endmatrix\right)^{-1})\right|=
\left|{rd-sc\over ad-bc}\right|
\leq{\|(r,s)\|\cdot\|(c,d)\|\over|ad-bc|}.$$
On the other hand let $(a_1,b_1)=(a,b)g$ and $(c_1,d_1)=(c,d)g$. Then $|a_1|\geq|b_1|$ and $|c_1|<|d_1|$ similarly as we noted at the end of Definition 3.2 therefore
$$\split|ad-bc|=&
\left|\det\left(\matrix a_1&b_1\\c_1&d_1\endmatrix\right)\cdot
\det(g^{-1})\right|=|a_1d_1-b_1c_1|\cdot
|\det(g)|^{-1}\\
=&|a_1d_1|\cdot|\det(g)|^{-1}\\
=&\|(a_1,b_1)\|\cdot\|(c_1,d_1)\|\cdot|\det(g)|^{-1}\\
\geq&\|(a,b)g\|\cdot\|(c,d)\|\cdot\|g^{-1}\|^{-1}
\cdot|\det(g)|^{-1}\text{.\ $\square$}\endsplit$$
\enddefinition
\definition{Definition 3.7} As a consequence of Proposition 3.5 the
function $E^G_{\goth m}(C,D,g,x,y)$ can be evaluated at $x=y=1$. Let
$$E^G_{\goth m}(C,D,g)=
E^G_{\goth m}(C,D,g,1,1)\in R[G]$$
for every $g\in GL_2(\Bbb A)$. In accordance with the previously introduced notation for every finite abelian group $G$ we let $I_G\triangleleft\Bbb Z[G]$ denote the augmentation ideal of $\Bbb Z[G]$ that is the kernel of the augmentation map $\Bbb Z[G]\rightarrow\Bbb Z$. There is an isomorphism $I_G/I_G^2=G$ induced by the map given by the rule $g\mapsto 1-g\in I_G$ for every $g\in G$.
\enddefinition
\proclaim{Proposition 3.8} Assume that $R=\Bbb Z$. Then we have $E^G_{\goth m}(C,D,g)\in I_G$ for every $g\in GL_2(\Bbb A)$.
\endproclaim
\definition{Proof} It will be sufficient to prove that $E^{\{1\}}_{\goth m}(C,D,g)=0$
where $\{1\}$ is the trivial group. We may assume again by bilinearity that $C=(\alpha,\beta)-(\gamma,\delta)$  and $D=(\epsilon,\iota)-(\kappa,\lambda)$ for some pairs $(\alpha,\beta)$, $(\gamma,\delta)$, $(\epsilon,\iota)$ and $(\kappa,\lambda)\in\Cal V_{\goth m}$. Pick again two elements $(r,s)\in U_{\goth m}(\alpha-\gamma,\beta-\delta,g)$ and $(u,v)\in V_{\goth m}(\epsilon-\kappa,\iota-\lambda,g)$. Then for every sufficiently large natural number $n$ we have:
$$\bigcup_{m\leq n}U_{\goth m}(\alpha,\beta,g,m)=\bigcup_{m\leq n}\{(a+r,b+s)|(a,b)\in
U_{\goth m}(\gamma,\delta,g,m)\}$$
and
$$\bigcup_{m\leq n}V_{\goth m}(\epsilon,\iota,g,m)=\bigcup_{m\leq n}\{(a+u,b+v)|(a,b)\in
V_{\goth m}(\kappa,\lambda,g,m)\}.$$
Hence we have:
$$|\bigcup_{m\leq n}U_{\goth m}(\alpha,\beta,g,m)|=|\bigcup_{m\leq n}
U_{\goth m}(\gamma,\delta,g,m)|,|U_{\goth m}(\alpha,\beta,g,n)|=|
U_{\goth m}(\gamma,\delta,g,n)|\text{\ and}$$
$$|\bigcup_{m\leq n}V_{\goth m}(\epsilon,\iota,g,m)|=|\bigcup_{m\leq n}
V_{\goth m}(\kappa,\lambda,g,m)|,|V_{\goth m}(\epsilon,\iota,g,n)|=|
V_{\goth m}(\kappa,\lambda,g,n)|$$
for every sufficiently large natural number $n$.  Therefore 
$$\split E^{\{1\}}_{\goth m}(C,D,g)=&
\sum_{m,n\in\Bbb Z}\bigg((|U_{\goth m}(\alpha,\beta,g,m)|-|
U_{\goth m}(\gamma,\delta,g,m)|)\\
&\cdot
(|V_{\goth m}(\epsilon,\iota,g,n)|-|V_{\goth m}(\kappa,\lambda,g,n)|)\bigg)\\
=&
\lim_{n\rightarrow\infty}\bigg(
(|\bigcup_{m\leq n}U_{\goth m}(\alpha,\beta,g,m)|-|\bigcup_{m\leq n}
U_{\goth m}(\gamma,\delta,g,m)|)\\
&\cdot(|\bigcup_{k\leq n}V_{\goth m}(\epsilon,\iota,g,k)|-
|\bigcup_{k\leq n}V_{\goth m}(\kappa,\lambda,g,k)|)\bigg)
\\=&0\text{.\ $\square$}\endsplit$$
\enddefinition
\definition{Definition 3.9} In accordance with the notation we introduced in Definition 1.8
let $\theta'\in G$ denote the class of any $\theta\in I_G$. For every $C$,
$D\in\Bbb Z[\Cal V_{\goth m}]_0$ and $N\in\Bbb Z$ let $\Cal E_{\goth m}(C,D,g,n)$ denote the
$F_{\infty}^*$-valued function:
$$\Cal E_{\goth m}(C,D,g,N)=
\prod\Sb m,n\leq n\\(\alpha,\beta)\in\Cal V_{\goth m}\\(\gamma,\delta)\in\Cal V_{\goth m}\endSb\prod\Sb(a,b)\in U_{\goth m}(\alpha,\beta,g,m)\\(c,d)\in
V_{\goth m}(\gamma,\delta,g,n)\endSb
\!\!\!\!\!\!\!\!\!
\det\left(\matrix a&b\\c&d\endmatrix\right)^{C(\alpha,\beta)D(\gamma,\delta)}.$$
Finally let $\Cal E_{\goth m}(C,D,g)$ denote the limit 
$$\Cal E_{\goth m}(C,D,g)=\lim_{N\rightarrow\infty}\Cal E_{\goth m}(C,D,g,N)$$
if the latter exists. The following claim is an immediate corollary to Lemma 3.6 and Proposition 3.8 using the same argument we used in the proof of Proposition 3.5.
\enddefinition
\proclaim{Proposition 3.10} The limit above exits and
$$E^G_{\goth m}(C,D,g)'=\Cal E_{\goth m}(C,D,g)^G\text{.\ $\square$}$$
\endproclaim

\head 4. The Kronecker limit formula
\endhead

\definition{Notation 4.1} We are going to use the notation we introduced in 1.2. For every connected rational subdomain $U$ of $\Bbb P^1$ the elements of $\partial U$ are called the boundary components of $U$, by slight abuse of language. Let $\Cal R(U)\subset\Cal O(U)$ denote the subalgebra of restrictions of rational functions holomorphic on $U$ and $\Cal R^*(U)$ denote
the group of invertible elements of this algebra. The group $\Cal
R^*(U)$ consists of rational functions which do not have poles or zeros
lying in $U$.
\enddefinition
\proclaim{Theorem 4.2} There is a unique map $\{\cdot,\cdot\}_D:\Cal
O^*(U)\times\Cal O^*(U)\rightarrow\Bbb C_{\infty}^*$ for every $D\in\partial U$,
called the rigid analytic regulator, with the following properties:

\noindent (i) For any two $f$, $g\in\Cal R^*(U)$ their regulator is:
$$\{f,g\}_D=\prod_{x\in D}\{f,g\}_x,$$

\noindent (ii) the regulator $\{\cdot,\cdot\}_D$ is bilinear
in both variables,

\noindent (iii) the regulator $\{\cdot,\cdot\}_D$ is alternating: 
$\{f,g\}_D\cdot\{g,f\}_D=1$,

\noindent (iv) if $f$, $1-f\in\Cal O(U)^*$, then $\{f,1-f\}_D$ is $1$,

\noindent (v) for each $f\in\Cal O_{\epsilon}(U)$ and $g\in\Cal O^*(U)$
we have $\{f,g\}_D\in \Cal U_{\epsilon}$.
\endproclaim
\definition{Proof} This is Theorem 2.2 of [22].\ $\square$ 
\enddefinition
\definition{Definition 4.3} If $U$ is still a connected rational subdomain
of $\Bbb P^1$, and $f$, $g$ are two meromorphic functions on $U$, then
 for all $x\in U$ the functions $f$ and $g$ have a power series expansion
around $x$, so in particular their tame symbol $\{f,g\}_x$ at $x$ is
defined. Let $\Cal M(U)$ denote the field of meromorphic functions of
$U$. The tame symbols extends to a homomorphism
$\{\cdot,\cdot\}_x:K_2(\Cal M(U))\rightarrow\Bbb C_{\infty}^*$. We define the
group $K_2(U)$ as the kernel of the direct sum of tame symbols:
$$\bigoplus_{x\in U}\{\cdot,\cdot\}_x:K_2(\Cal M(U))\rightarrow
\bigoplus_{x\in U}\Bbb C_{\infty}^*.$$
Let $k=\sum_if_i\otimes g_i\in K_2(U)$, where $f_i$, $g_i\in\Cal M(U)$,
and let $D\in\partial U$. Let moreover $Y$ be a connected rational
subdomain of $U$ such that $f_i$, $g_i\in\Cal O^*(Y)$ for all $i$ and
$\partial U\subseteq\partial Y$. Define the rigid analytical regulator
$\{k\}_D$ by the formula:
$$\{k\}_D=\prod_i\{f_i|_Y,g_i|_Y\}_D.$$
\enddefinition
Theorem 1.3 is based on the previous result and the following theorem:
\proclaim{Theorem 4.4} $(i)$ For each $k\in K_2(U)$ the rigid analytical
regulator $\{k\}_D$ is well-defined, and it is a homomorphism
$\{\cdot\}_D:K_2(U)\rightarrow\Bbb C_{\infty}^*$,

\noindent $(ii)$ for any two functions $f$, $g\in\Cal O^*(U)$ we have
$\{f\otimes g\}_D=\{f,g\}_D$,

\noindent $(iii)$ for every $k\in K_2(U)$ the product of all regulators
on the boundary components of $U$ is equal to 1:
$$\prod_{D\in\partial U}\{k\}_D=1,$$

\noindent $(iv)$ for every connected subdomain $Y\subseteq U$, boundary component
$D\in\partial Y\cap\partial U$ and $k\in K_2(\Cal M(U))$ we have:
$$\{k|_Y\}_D=\{k\}_D.$$
\endproclaim
\definition{Proof} This is Theorem 3.2 of [22].\ $\square$ 
\enddefinition
\definition{Definition 4.5} For every $\rho\in GL_2(F_{\infty})$ and
$z\in\Bbb P^1$ let $\rho(z)$ denote the image of $z$ under the M\"obius transformation corresponding to $\rho$. Let moreover $D(\rho)$ denote the open disc
$$D(\rho)=\{z\in\Bbb P^1(\Bbb C_{\infty})|1<|\rho^{-1}(z)|\}.$$
Let $\Cal D$ denote the set of open discs of the form $D(\rho)$ where $\rho\in GL_2(F_{\infty})$. For each $D\in\Cal D$ let $D(F_{\infty})$ denote  $D\cap\Bbb P^1(F_{\infty})$. Let $\Cal P$ denote those subsets $S$ of
$\Cal D$ such that the sets $D(F_{\infty})$, $D\in S$ form a pair-wise
disjoint partition of $\Bbb P^1(F_{\infty})$. For each $S\in\Cal P$ let
$\Omega(S)$ denote the unique connected rational subdomain defined over $F_{\infty}$ with the property $\partial\Omega(S)=S$. Let $\Omega$ denote the rigid analytic upper half plane, or Drinfeld's upper half plane over $F_{\infty}$. The set of points of $\Omega$ is $\Bbb C_{\infty}-F_{\infty}$, denoted also by $\Omega$ by abuse of notation. Recall that a function $f:\Omega\rightarrow\Bbb C_{\infty}$ is holomorphic if the restriction of $f$ onto $\Omega(S)$ is holomorphic for every $S\in\Cal P$. Let $\Cal O(\Omega)$ and $\Cal M(\Omega)$ denote the $\Bbb C_{\infty}$-algebra of holomorphic functions and the field of meromorphic function of $\Omega$, respectively. The latter is of course the quotient field of the former.  We define $K_2(\Omega)$ as the intersection of the kernels of all the tame symbols $\{\cdot,\cdot\}_x$ inside $K_2(\Cal M(\Omega))$ where $x$ runs through the set $\Omega$.  By part $(iv)$ of Theorem 4.4 for each $k\in K_2(\Omega)$ the value  $\{k\}(\rho)=\{k|_{
 \Omega(S)}\}_{D(\rho)}$, where $\rho\in GL_2(F_{\infty})$ and $D(\rho)\in S\in\Cal P$, is independent of the choice of $S$. We define the regulator $\{k\}:GL_2(F_{\infty})\rightarrow\Bbb C_{\infty}^*$ of $k$ as the function given by this rule.
\enddefinition
\proclaim{Lemma 4.6} Let $\rho=\left(\smallmatrix x&y\\0&1\endsmallmatrix\right)$ where $x\in F^*_{\infty}$ and $y\in F_{\infty}$. Then for every $0\neq(a,b)\in F^2_{\infty}$ and $0\neq(c,d)\in F^2_{\infty}$ the following holds:
\roster
\item"$(i)$" if $(a,b)\rho\in F^2_<$ and $(c,d)\rho\in F^2_<$ then
$$\{(az+b)\otimes(cz+d)\}_{D(\rho)}=1,$$
\item"$(ii)$" if $(a,b)\rho\notin F^2_<$ and $(c,d)\rho\notin F^2_<$ then $a\neq 0$,  $b\neq 0$ and
$$\{(az+b)\otimes(cz+d)\}_{D(\rho)}=b/a,$$
\item"$(iii)$" if $(a,b)\rho\notin F^2_<$ and $(c,d)\rho\in F^2_<$ then $a\neq 0$ and 
$$\{(az+b)\otimes(cz+d)\}_{D(\rho)}={1\over a}
\det\left(\matrix a&b\\c&d\endmatrix\right),$$
\item"$(iv)$" if $(a,b)\rho\in F^2_<$ and $(c,d)\rho\notin F^2_<$ then $b\neq 0$ and 
$$\{(az+b)\otimes(cz+d)\}_{D(\rho)}=c
\det\left(\matrix a&b\\c&d\endmatrix\right)^{-1}.$$
\endroster
\endproclaim
\definition{Proof} Let $D(\rho)^c$ denote the complement of $D(\rho)$ in $\Bbb P^1$. Obviously
$$D(\rho)^c=\{z\in\Bbb C_{\infty}||z-y|\leq|x|\}.$$
Hence $(a,b)\rho\in F^2_<$ if and only if the polynomial $az+b$ has no zeros in $D(\rho)^c$ and
$(a,b)\rho\notin F^2_<$ if and only if $a\neq0$ and the polynomial $az+b$ does have a zero in
$D(\rho)^c$. By Weil's reciprocity law:
$$\{(az+b)\otimes(cz+d)\}_{D(\rho)}^{-1}=\prod_{t\in D(\rho)^{-1}}\{az+b,cz+d\}_t.$$
Since the tame symbol of $az+b$ and $cz+d$ at $t\in\Bbb C_{\infty}$ is $1$ if neither $az+b$ nor $cz+d$ has a zero at $t$ claim $(i)$ is clear. In the second case both $az+b$ and $cz+d$ has a single pole at $\infty$ but has no zero in $D(\rho)$ therefore
$$\{(az+b)\otimes(cz+d)\}_{D(\rho)}=\{az+b,cz+d\}_{\infty}=b/a.$$
Claim $(iv)$ follows from claim $(iii)$ by the antisymmetry of the regulator. In the latter case $az+b$ has a single zero in $D(\rho)^c$ and $cz+d$ has no zeros in $D(\rho)^c$ hence
$$\{(az+b)\otimes(cz+d)\}_{D(\rho)}=\{az+b,cz+d\}_{-b/a}^{-1}={1\over a}
\det\left(\matrix a&b\\c&d\endmatrix\right)\text{.\ $\square$ }$$
\enddefinition
\definition{Definition 4.7} We are going to need a mild extension of the regulator we introduced in Definition 4.5. Let $K_2(GL_2(\Bbb A_f)\times\Omega)$ denote the set of functions $k:GL_2(\Bbb A_f)\rightarrow K_2(\Omega)$. We define the regulator of an element $k\in K_2(GL_2(\Bbb A_f)\times\Omega)$ as the function $\{k\}:GL_2(\Bbb A)
\rightarrow\Bbb C^*_{\infty}$ given by the rule $\{k\}(g)=\{k(g_f)\}(g_{\infty})$ for every $g\in GL_2(\Bbb A)$. Since the set
$K_2(GL_2(\Bbb A_f)\times\Omega)$ consists of functions taking values in the group $K_2(\Omega)$ it is equipped with a group structure
whose operation will be denoted by addition. Let $\Cal O^*(GL_2(\Bbb A_f)\times\Omega)$ denote the set of functions $u:GL_2(\Bbb A_f)\times\Omega
\rightarrow\Bbb C^*_{\infty}$ which are holomorphic in the second variable. Then there is a bilinear map:
$$\otimes:\Cal O^*(GL_2(\Bbb A_f)\times\Omega)\times\Cal O^*(GL_2(\Bbb A_f)\times\Omega)
\rightarrow K_2(GL_2(\Bbb A_f)\times\Omega)$$
given by the rule $(u\otimes v)(g)=u|_{g\times\{\cdot\}}\otimes v|_{g\times\{\cdot\}}$ for every $g\in GL_2(\Bbb A_f)$. For each $(\alpha,\beta)\in(\Cal O_f/\goth m)^2$, and $N$ positive integer let $\epsilon_{\goth m}(\alpha,\beta,N)(g,z)$ denote the function:
$$\epsilon_{\goth m}(\alpha,\beta,N)(g,z)=
\prod_{n\leq N}\left(
\prod_{(a,b)\in W_{\goth m}(\alpha,\beta,g,n)}
\!\!\!\!\!\!\!(az+b)\cdot\!\!\!\!\!\!\!\!\!\!\!
\prod_{(c,d)\in W_{\goth m}(0,0,g,n)}\!\!\!\!\!\!\!
(cz+d)^{-1}\right).$$
on the set $GL_2(\Bbb A_f)\times\Omega$. The latter is clearly holomorphic in the second variable. 
\enddefinition
\proclaim{Lemma 4.8}  The limit
$$\epsilon_{\goth m}(\alpha,\beta)(g,z)=\lim_{N\rightarrow\infty}
\epsilon_{\goth m}(\alpha,\beta,N)(g,z)$$
converges uniformly in $z$ on every admissible open subdomain of $\Omega$
for every fixed $g$ and defines a function ho\-lo\-mor\-phic in the second
variable.
\endproclaim
\definition{Proof} See Lemma 4.5 of [20] on pages 145-146.\ $\square$
\enddefinition
\definition{Definition 4.9} For every $C\in\Bbb Z[\Cal V_{\goth m}]_0$ let $\epsilon_{\goth m}(C,g,z)$ denote the function:
 $$\prod_{(\alpha,\beta)\in\Cal V_{\goth m}}
 \epsilon_{\goth m}(\alpha,\beta)(g,z)^{C(\alpha,\beta)}$$
 on the set $GL_2(\Bbb A_f)\times\Omega$. For every $C$, $D\in\Bbb Z[\Cal V_{\goth m}]_0$ let $\kappa_{\goth m}(C,D)$ denote the element:
$$\epsilon_{\goth m}(C,g,z)\otimes\epsilon_{\goth m}(D,g,z)$$
of the set $K_2(GL_2(\Bbb A_f)\times\Omega)$.
\enddefinition
\proclaim{Kronecker Limit Formula 4.10} For all $g\in GL_2(\Bbb A)$ we have:
$$\{\kappa_{\goth m}(C,D)\}(g)^G
=(E^G_{\goth m}(C,D,g)-E^G_{\goth m}(D,C,g))'.$$
\endproclaim
\definition{Proof} Assume first that $g_{\infty}=\left(\smallmatrix x&y\\0&1\endsmallmatrix
\right)$ for some $x\in F^*_{\infty}$ and $y\in F_{\infty}$. By Proposition 3.10 it will be sufficient to prove that
$$\{\prod_{(\alpha,\beta)\in\Cal V_{\goth m}}
\epsilon_{\goth m}(\alpha,\beta,N)(g,z)^{C(\alpha,\beta)}\otimes
\prod_{(\gamma,\delta)\in\Cal V_{\goth m}}
\epsilon_{\goth m}(\gamma,\delta,N)(g,z)^{D(\gamma,\delta)}\}_{D(g_{\infty})}$$
is equal to
$$\Cal E_{\goth m}(C,D,g,N)\cdot \Cal E_{\goth m}(D,C,g,N)^{-1}$$
for every sufficiently large $N$. By bilinearity and Lemma 4.6 the regulator in the left hand side of the equation that we wish to prove is equal to
$$\split=\prod\Sb m,n\leq N\\(\alpha,\beta)\in\Cal V_{\goth m}\\(\gamma,\delta)\in
\Cal V_{\goth m}\endSb
\bigg(&
\prod\Sb(a_1,b_1)\in U_{\goth m}(\alpha,\beta,g,m)\\(c_1,d_1)\in
V_{\goth m}(\gamma,\delta,g,n)\endSb
\!\!\!\!\!\!\!\!\!
a_1^{-1}\det\left(\matrix a_1&b_1\\c_1&d_1\endmatrix\right)\\
&\cdot
\prod\Sb(a_2,b_2)\in V_{\goth m}(\alpha,\beta,g,m)\\(c_2,d_2)\in
U_{\goth m}(\gamma,\delta,g,n)\endSb
\!\!\!\!\!\!\!\!\!
c_2\det\left(\matrix a_2&b_2\\c_2&d_2\endmatrix\right)^{-1}\\
&\cdot\prod\Sb(a_3,b_3)\in U_{\goth m}(\alpha,\beta,g,m)\\(c_3,d_3)\in
U_{\goth m}(\gamma,\delta,g,n)\endSb
\!\!\!\!\!\!\!\!\!a_3^{-1}c_3\bigg)^{C(\alpha,\beta)D(\gamma,\delta)}
.\endsplit$$
Therefore what we need to show is:
$$\prod\Sb m,n\leq N\\(\alpha,\beta)\in\Cal V_{\goth m}\\(\gamma,\delta)\in
\Cal V_{\goth m}\endSb\bigg(
\prod_{(a,b)\in U_{\goth m}(\alpha,\beta,g,m)}
\!\!\!\!\!\!\!\!\!a^{-|W_{\goth m}(\gamma,\delta,g,n)|}
\cdot\!\!\!\!\!\!\!\!\!
\prod_{(c,d)\in U_{\goth m}(\gamma,\delta,g,n)}
\!\!\!\!\!\!\!\!\!c^{|W_{\goth m}(\alpha,\beta,g,m)|}
\bigg)^{C(\alpha,\beta)D(\gamma,\delta)}=1.$$
The latter follows from the fact that for every $C\in\Bbb Z[\Cal V_{\goth m}]_0$ and for every sufficiently large $N$ the equation:
$$\sum\Sb n\leq N\\(\alpha,\beta)\in\Cal V_{\goth m}\endSb
C(\alpha,\beta)|W_{\goth m}(\alpha,\beta,g,n)|=0$$
holds. On the other hand the latter has been already shown in the course of the proof of Proposition 3.8 (at least in the special case when
$C=(\alpha,\beta)-(\gamma,\delta)$ for some $(\alpha,\beta)$,
$(\gamma,\delta)\in\Cal V_{\goth m}$ but the general case follows at
once from this one by linearity). Let us consider now the general case. First note that both sides of the equation in the theorem above are right-invariant with respect to $Z(F_{\infty})\Gamma_{\infty}$ hence if the claim is true for $g$ then it is true for $gz$ as well for every $z\in Z(F_{\infty})\Gamma_{\infty}$. Let $\Pi\in GL_2(F_\infty)$ be the matrix whose diagonal entries are zero, and its lower left and upper right entry is $\pi$ and $1$,
respectively, where $\pi$ is a uniformizer of $F_{\infty}$. Then for every
$\rho\in GL_2(F_{\infty})$ the open disks $D(\rho)$ and $D(\rho\Pi)$ are complementary in $\Bbb P^1$ hence claim $(iii)$ of Theorem 4.4 (Weil's reciprocity law) implies that
$$\{\kappa_{\goth m}(C,D)\}(g)^G\cdot\{\kappa_{\goth m}(C,D)\}(g\Pi)^G=1$$
for every $g\in GL_2(\Bbb A)$. A matrix $\rho\in GL_2(F_{\infty})$ can be written as a product $\rho=\left(\smallmatrix x&y\\0&1\endsmallmatrix\right)z$ where $x\in F^*_{\infty}$,
$y\in F_{\infty}$ and $z\in Z(F_{\infty})\Gamma_{\infty}$ if and only if $\infty\in D(\rho)$.
Since either $D(\rho)$ or $D(\rho\Pi)$ contains the point $\infty$ it will be sufficient to prove that the identity above holds for
$E^G_{\goth m}(C,D,g)'/E^G_{\goth m}(D,C,g)'$ as well.
But
$$U_{\goth m}(\alpha,\beta,g\Pi)=V_{\goth m}(\alpha,\beta,g)
\text{\ and\ }
V_{\goth m}(\alpha,\beta,g\Pi)=U_{\goth m}(\alpha,\beta,g),$$
therefore $E^G_{\goth m}(C,D,g)=(-1)^GE^G_{\goth m}(D,C,g\Pi)$
for any $g\in GL_2(\Bbb A)$, so the latter is obvious.\ $\square$
\enddefinition

\heading 5. Modular units and Hecke operators
\endheading

\definition{Notation 5.1} Let $A=\Cal O_f\cap F$: it is a Dedekind
domain. The ideals of $A$ and the effective divisors on $X$ with support
away from $\infty$ are in a bijective correspondence. These two sets will
be identified in all that follows. For any non-zero ideal $\goth m\triangleleft A$ let
$Y(\goth m)$ denote the coarse moduli for parameterizing Drinfeld modules of rank two over $A$
of general characteristic with full level $\goth m$-structure. It is an affine algebraic curve
defined over $F$. For every Drinfeld module $\phi:A\rightarrow C\{\tau\}$  of rank two
equipped with a full level $\goth m$-structure $\iota:(A/\goth m)^2\rightarrow C$, where $C$
is an $F$-algebra, let $u_{\phi,\iota}:Y(\goth m)\rightarrow\text{Spec}(C)$ be the universal
map. 
\enddefinition
\proclaim{Lemma 5.2} There is a unique element $\epsilon_{\goth m}(D)
\in\Gamma(Y(\goth m),\Cal O^*)$ for every $D\in\Bbb Z[\Cal V_{\goth m}]_0$ such that
$$u^*_{\phi,\iota}(\epsilon_{\goth m}(D))
=\prod_{(\alpha,\beta)\in\Cal V_{\goth m}}\iota(\alpha,\beta)^{D(\alpha,\beta)}\in C$$
for every $C$, $\phi$ and $\iota$ as above. 
\endproclaim
\definition{Proof}  We may assume that $\goth m$ is a proper ideal without the loss of generality.
For any non-zero ideal $\goth m\triangleleft A$ let $H(\goth m)$ denote
$\Gamma(Y(\goth m),\Cal O)$. We may assume by linearity that
$D=(\alpha,\beta)-(\gamma,\delta)$ for some $(\alpha,\beta)$,
$(\gamma,\delta)\in\Cal V_{\goth m}$. Let $(\phi,\iota)$ and $(\psi,\kappa)$ be ordered pairs
of two Drinfeld
modules $\phi$ and $\psi$ of rank two over $C$ equipped with a full
level $\goth m$-structure $\iota$ and $\kappa$, respectively. Recall that
$(\phi,\iota)$ and $(\psi,\kappa)$ are isomorphic if there is an
isomorphism $j:\Bbb G_a\rightarrow\Bbb G_a$ between $\phi$ and $\psi$
such that the composition $j\circ\iota$ is equal to $\kappa$. As $j$ is just a multiplication
by a scalar we get that the element $\iota(\alpha,\beta)/\iota(\gamma,\delta)$ depends only on
the isomorphism class of the pair $(\phi,\iota)$. In particular the claim is
obvious when the moduli scheme $Y(\goth m)$ is fine because we have
$\epsilon_{\goth m}(D)=\iota_{\goth m}(\alpha,\beta))/\iota_{\goth m}(\gamma,\delta)$
in this case where the map
$\iota_{\goth m}:(A/\goth m)^2\rightarrow H(\goth m)$ is the
universal  full level $\goth m$-structure for the universal Drinfeld
module $\phi_{\goth m}:A\rightarrow H(\goth m)\{\tau\}$ over $Y(\goth m)$. The latter
holds if $\goth m$ has at least two prime divisors. In general the
universal map
$H(\goth m)\rightarrow\bigoplus_{\goth p\nmid\goth m}H(\goth m\goth p)$
is an \'etale injection, so it is faithfully flat. Therefore the sequence
$$0\rightarrow H(\goth m)\rightarrow\bigoplus_{\goth p\nmid\goth m}
H(\goth m\goth p)\rightrightarrows
\left(\bigoplus_{\goth p\nmid\goth m}H(\goth m\goth p)\right)
\otimes_{H(\goth m)}
\left(\bigoplus_{\goth p\nmid\goth m}H(\goth m\goth p)\right)$$
is exact by Proposition 2.18 of [18], pages 16-17. For every
$(\kappa,\lambda)\in\Cal V_{\goth m}$ and prime ideal $\goth p\nmid\goth m$ let
$(\kappa,\lambda,\goth p)$ denote the unique element of $\Cal V_{\goth{mp}}$ such that
$(\kappa,\lambda,\goth p)\equiv(\kappa,\lambda)\mod\goth m$ and
$(\kappa,\lambda,\goth p)\equiv(0,0)\!\!\mod\goth p$. Moreover for every prime ideal
$\goth p\nmid\goth m$ let $D(\goth p)\in\Bbb Z[\Cal V_{\goth{mp}}]_0$ denote the element
$(\alpha,\beta,\goth p)-(\gamma,\beta,\goth p)$. Then the element
$$\bigoplus_{\goth p\nmid\goth m}
\epsilon_{\goth{mp}}(D(\goth p))
\in\bigoplus_{\goth p\nmid\goth m} H(\goth{mp})$$
is in the kernel of the second map in the exact sequence above, therefore
it is the image of a unique element $\epsilon_{\goth m}(D)\in H(\goth m)$ which satisfies the
required property.\ $\square$
\enddefinition
\definition{Definition 5.3} The group
$GL_2(F)$ acts on the product $GL_2(\Bbb A_f)\times
\Omega$ on the left by acting on the first factor via the natural
embedding and on Drinfeld's upper half plane via M\"obius
transformations. The group $\Bbb K_f(\goth m)=\Bbb K(\goth m)\cap
GL_2(\Cal O_f)$ acts on the right of this product by acting on the first
factor via the regular action. Since the quotient set $GL_2(F)\backslash
GL_2(\Bbb A_f)/\Bbb K_f(\goth m)$ is finite, the set
$$GL_2(F)\backslash GL_2(\Bbb A_f)\times\Omega/\Bbb K_f(\goth m)$$
is the disjoint union of finitely many sets of the form $\Gamma\backslash
\Omega$, where $\Gamma$ is a subgroup of $GL_2(F)$ of the form
$GL_2(F)\cap g\Bbb K_f(\goth m)g^{-1}$ for some $g\in GL_2(\Bbb A_f)$. As
these groups act on $\Omega$ discretely, the set above naturally
has the structure of a rigid analytic curve. Let $Y(\goth m)_{F_{\infty}}$ also denote the underlying rigid analytical space of the base change of $Y(\goth m)$ to $F_{\infty}$ by abuse of notation.
\enddefinition
\proclaim{Theorem 5.4} There is a rigid-analytical isomorphism:
$$Y(\goth m)_{F_{\infty}}
\cong GL_2(F)\backslash GL_2(\Bbb A_f)\times\Omega/
\Bbb K_f(\goth m).$$
\endproclaim
\definition{Proof} See [4], Theorem 6.6.\ $\square$
\enddefinition
\proclaim{Proposition 5.5} For every $D\in\Bbb Z[\Cal V_{\goth m}]_0$ the function
corresponding to $\epsilon_{\goth m}(D)$ under the isomorphism of Theorem 5.4 above is the
function $\epsilon_{\goth m}(D,g,z)$ introduced in Definition 4.9.
\endproclaim
\definition{Proof} First we are going to recall the map underlying the isomorphism of Theorem
5.4 on $\Bbb C_{\infty}$-valued points. For every $(g,z)\in GL_2(\Bbb A_f)\times\Omega$ let
$e_{(g,z)}(w)$ denote the corresponding exponential function:
$$e_{(g,z)}(w)=z\prod_{(a,b)\in W_{\goth m}(0,0,g_f)}(1-{w\over az+b}
).$$
The infinite product above is converging absolutely and defines an entire function. The exponential $e_{(g,z)}$ uniformizes a Drinfeld module $\phi_{(g,z)}$ over $\Bbb C_{\infty}$ which is equipped with a full level $\goth m$-structure
$\iota$ given by the formula:
$$\iota(\alpha,\beta)=e_{(g,z)}(az+b)\text{\ where }(a,b)\in W_{\goth m}(\alpha,\beta,g_f)$$
for every $(\alpha,\beta)\in(\Cal O_f/\goth m)^2$ independent of the choice of $(a,b)$. Since
obviously we have $\iota(\alpha,\beta)=\epsilon_{\goth m}(\alpha,\beta,g,z)$ the claim is now
clear.\ $\square$
\enddefinition
\definition{Definition 5.6} Let $M$ be an abelian group and let $\goth n$ be any effective
divisor on $X$. By an $M$-valued
automorphic form over $F$ of level $\goth n$ (and trivial central character) we mean a
locally constant function $\phi:GL_2(\Bbb A)\rightarrow M$ satisfying $\phi(\gamma gkz)=
\phi(g)$ for all $\gamma\in GL_2(F)$, $z\in Z(\Bbb A)$, and $k\in\Bbb K_0(\goth n)$,
where 
$$\Bbb K_0(\goth n)=\{\left(\matrix a&b\\c&d\endmatrix\right)
\in GL_2(\Cal O)|c\equiv0\text{\ mod }\goth n\}.$$
Let $\Cal A(\goth n,M)$ denote the $\Bbb Z$-module of $M$-valued automorphic forms
of level $\goth n$. Now let $\goth n$
be an effective divisor on $X$ whose support does not contain $\infty$.
Let $\Cal H(\goth n,M)$ denote the $\Bbb Z$-module of
automorphic forms $f$ of level $\goth n\infty$  satisfying the following two identities:
$$\phi(g\left(\matrix 0&1\\\pi&0\endmatrix\right))
=-\phi(g),(\forall g\in GL_2(\Bbb A)),$$
and
$$\phi(g\left(\matrix 0&1\\1&0\endmatrix\right))+
\sum_{\epsilon\in\bold f_{\infty}}
\phi(g\left(\matrix1&0\\\epsilon&1\endmatrix\right))=0,
(\forall g\in GL_2(\Bbb A)),$$
where $\pi$ is a uniformizer in $F_{\infty}$ and we consider
$GL_2(F_{\infty})$ as a subgroup of $GL_2(\Bbb A)$ and we understand the
product of their elements accordingly. Such automorphic forms are called
harmonic.
\enddefinition
\definition{Definition 5.7} Let $\goth m$, $\goth n$ be effective divisors
of $X$. Define the set:
$$H(\goth m,\goth n)=
\{\left(\matrix a&b\\c&d\endmatrix\right)\in GL_2(\Bbb A)|a,b,c,d\in\Cal
O,(ad-cb)=\goth m,
\goth n\supseteq(c),(d)+\goth n=\Cal O\}.$$
The set $H(\goth m,\goth n)$ is compact and it is a  double $\Bbb K_0
(\goth n)$-coset, so it is a disjoint union of finitely many right
$\Bbb K_0 (\goth n)$-cosets. Let $R(\goth m,\goth n)$ be a set of
representatives of these cosets. For any $\phi\in\Cal A(\goth n,R)$
define the function $T_{\goth m}(\phi)$ by the formula: 
$$T_{\goth m}(\phi)(g)=\sum_{h\in R(\goth m,\goth n)}\phi(gh).$$
It is easy to check that $T_{\goth m}(\phi)$ is independent of the
choice of $R(\goth m,\goth n)$ and $T_{\goth m}(\phi)\in\Cal A(\goth
n,M)$ as well. So we have an $\Bbb Z$-linear operator $T_{\goth m}:\Cal
A(\goth n,M)\rightarrow\Cal A(\goth n,M)$. 
\enddefinition
\definition{Definition 5.8} Let $\goth X$ be a Hausdorff topological
space. For any $R$ commutative ring let $\Cal M(\goth X,R)$ denote the set
of $R$-valued finitely additive measures on the open and compact
subsets of $\goth X$.  For every $M$ abelian group let $\Cal C_0(\goth
X,M)$ denote the group of compactly supported locally constant functions
$f:\goth X\rightarrow M$. For every $f\in\Cal C_0 (\goth X,M)$ and $\mu\in\Cal
M(\goth X,R)$ we define the modulus $\mu(f)$ of $f$ with respect to $\mu$
as the $\Bbb Z$-submodule of $R$ generated by the elements
$\mu(f^{-1}(g))$, where $0\neq g\in M$. We also define the integral of $f$ on
$\goth X$ with respect to $\mu$ as the sum:
$$\int_{\goth X}f(x)d\mu(x)=\sum_{g\in M}g\otimes\mu(f^{-1}(g))\in
M\otimes\mu(f).$$
Of course this definition is nothing more than a convenient
formalism. For its elementary properties see Lemma 5.2 of [22]. Let $M$ be
a $\Bbb Q$-vector space and let $\phi$ be an element of $\Cal A(\goth n,M)$.
If for all $g\in GL_2(\Bbb A)$:
$$\int_{F\backslash\Bbb A}\phi(\left(\matrix 1&x\\0&1\endmatrix
\right)g)d\mu(x)=0,$$
where $\mu$ is the normalized Haar measure on $F\backslash\Bbb A$
such that $\mu(F\backslash\Bbb A)=1$ we call $\phi$ a cusp form. Let
$\Cal A_0(\goth n,M)$ (respectively $\Cal H_0(\goth n,M)$) denote the
$\Bbb Q$-module of $M$-valued cuspidal automorphic forms
(respectively cuspidal harmonic forms) of level $\goth n$
(resp. of level $\goth n\infty$).
\enddefinition
\definition{Notation 5.9} For any ideal $\goth n\triangleleft A$ let
$Y_0(\goth n)$ denote the coarse moduli scheme for rank two Drinfeld
modules of general characteristic equipped with a Hecke level-$\goth n$ structure. It is an affine algebraic curve defined over $F$.
Let $X_0(\goth n)$ denote the unique irreducible smooth
projective curve over $F$ which contains $Y_0(\goth n)$ as
an open subvariety. For every proper ideal $\goth m\triangleleft A$ there is a $\goth m$-th Hecke correspondence on the Drinfeld modular curve $X_0(\goth n)$ which in turn induces an endomorphism of the Jacobian $J_0(\goth n)$ of the curve $X_0(\goth n)$, called
the Hecke operator $T_{\goth m}$ (for a detailed description see for example [6] or [7].) The $\goth m$-th Hecke
correspondence also induces a pair of compatible homomorphisms:
$$T_{\goth m}:H^2_{\Cal M}(X_0(\goth n)_L,K(2))\rightarrow 
H^2_{\Cal M}(X_0(\goth n)_L,K(2))$$
and
$$T_{\goth m}:H^2_{\Cal M}(Y_0(\goth n)_L,K(2))\rightarrow 
H^2_{\Cal M}(Y_0(\goth n)_L,K(2))$$
for every number field $K$ and for every $L\supseteq F$ extension. These operators are denoted by the same symbol we use for the operators introduced in Definition 5.7, but this will not cause confusion as we will see. For the moment it is sufficient to remark that they act on different objects.
\enddefinition
\definition{Notation 5.10} Let $\pi(\goth n):Y(\goth n)\rightarrow Y_0(\goth n)$ be the map induced by the forgetful functor
which assigns to every Drinfeld module $\phi:A\rightarrow C\{\tau\}$ of rank two equipped with a full level $\goth m$-structure $\iota:(A/\goth m)^2\rightarrow C$, where $C$ is an $F$-algebra, the Drinfeld module $\phi$ equipped with the Hecke level-$\goth n$
structure generated by $\iota(0,1)$. Hence $Y_0(\goth n)$ also has a rigid analytic uniformization of the kind described in Theorem 5.4 where the role of the group $\Bbb K(\goth n)$ is played by $\Bbb K_0(\goth n)$. Hence we may evaluate the regulator introduced in Definition 4.7 on the pull-back of the elements of 
$H^2_{\Cal M}(Y_0(\goth n)_{F_{\infty}},\Bbb Z(2))$ with respect to this uniformization. Let $\{\cdot\}$ denote also the unique $K$-linear extension to $H^2_{\Cal M}(Y_0(\goth n)_{F_{\infty}},K(2))$ of this regulator for every number field $K$ by abuse of notation.
\enddefinition
For the rest of the paper we assume that
$F=\Bbb F_q(T)$ is the rational function field of transcendence degree one over $\Bbb F_q$,
where $T$ is an indeterminate, and $\infty$ is the point at infinity on
$X=\Bbb P^1_{\Bbb F_q}$. 
\proclaim{Proposition 5.11} For every $k\in H^2_{\Cal M}
(Y_0(\goth n)_{F_{\infty}},K(2))$ following holds:
\roster
\item"$(i)$" we have $\{k\}\in\Cal H(\goth n,F_{\infty}^*\otimes K)$,
\item"$(ii)$" we have $\{k\}\in\Cal H_0(\goth n,F_{\infty}^*\otimes K)$ when $k\in H^2_{\Cal M}(X_0(\goth n)_{F_{\infty}},K(2))$,
\item"$(iii)$" we have $\{T_{\goth m}(k)\}=T_{\goth m}\{k\}$ for every $\goth m
\triangleleft A$ proper ideal.
\endroster
\endproclaim
\definition{Proof} By definition and the invariance theorem of [22] the
regulator $\{k\}$ is left $GL_2(F)$-invariant and right $\Bbb K_0(\goth n\infty)Z(F_{\infty})$-invariant. By our assumptions on $F$ and $\infty$ we have $F^*\Cal O_f^*=\Bbb A_f^*$ hence $\{k\}$ is also $Z(\Bbb A)$-invariant. Therefore it is an element of $\Cal A(\goth n\infty,F_{\infty}^*\otimes K)$. By claim $(iii)$ of Theorem 4.4 the additional conditions of Definition 5.6 also hold for $\{k\}$ as the following two sets of disks:
$$D(\rho),D(\rho\left(\matrix 0&1\\\pi&0\endmatrix\right))
\text{\ and\  }D(\rho\left(\matrix 0&1\\1&0\endmatrix\right)),
\{D(\rho\left(\matrix1&0\\\epsilon&1\endmatrix\right))|\epsilon\in\bold f_{\infty}\}$$
give a pair-wise disjoint covering of the set $\Bbb P^1(F_{\infty})$ for every
$\rho\in GL_2(F_{\infty})$. Claim $(i)$ is proved. The second claim is an immediate
consequence of Theorem 6.3 of [22]. Finally let us concern ourselves with the proof
of claim $(iii)$. For every $h\in GL_2(\Bbb A_f)$ let
$h:GL_2(\Bbb A_f)\times\Omega\rightarrow GL_2(\Bbb A_f)\times\Omega$ simply denote the
map given by the rule $(g,z)\mapsto (gh,z)$ for every $g\in GL_2(\Bbb A_f)$ and $z\in\Omega$. By slight abuse of notation let the same symbol  denote unique the map $h:Y(\goth n)_{F_{\infty}}\rightarrow Y(\goth n)_{F_{\infty}}$ which satisfies the
relation $h\circ\pi(\goth n)=\pi(\goth n)\circ h$. Let $R(\goth m,\goth n)
\subset H(\goth m,\goth n)$ be a set of representatives which is also a subset
of $GL_2(\Bbb A_f)$. Then we have:
$$\pi(\goth n)^*(T_{\goth m}(k))=\sum_{h\in R(\goth m,\goth n)}h^*(k)$$
in $H^2_{\Cal M}(Y(\goth n)_{F_{\infty}},K(2))$. Hence by the invariance theorem (Theorem 3.11 of [22]) for every $g\in GL_2(\Bbb A)$ we have:
$$\{T_{\goth m}(k)\}(g)=\sum_{h\in R(\goth m,\goth n)}\{h^*(k)\}(g)=
\sum_{h\in R(\goth m,\goth n)}\{k\}(gh)=T_{\goth m}\{k\}(g)$$
as we claimed.\ $\square$
\enddefinition
Let $L\subset F_{\infty}$ be a finite extension of $F$ and let
$$\{\cdot\}:H^2_{\Cal M}(X_0(\goth n)_L,K(2))
\rightarrow F^*_{\infty}\otimes K$$
denote also the composition of the homomorphism
$$H^2_{\Cal M}(X_0(\goth n)_L,K(2))\rightarrow H^2_{\Cal M}(X_0(\goth n)_{F_{\infty}},K(2))$$
induced by the functoriality of motivic cohomology and the homomorphism
$\{\cdot\}$. Let $C\subset X_0(\goth n)\times X_0(\goth n)$ be a correspondence and let 
$$C_*:H^2_{\Cal M}(X_0(\goth n)_L,K(2))\rightarrow 
H^2_{\Cal M}(X_0(\goth n)_L,K(2))$$
denote the homomorphism induced by $C$.
\proclaim{Lemma 5.12} We have $\{C_*(k)\}=0$ for every
$k\in H^2_{\Cal M}(X_0(\goth n)_L,K(2))$ if the endomorphism
$J(C):J_0(\goth n)\rightarrow J_0(\goth n)$ induced by $C$ is zero.
\endproclaim
\definition{Proof} We may assume without the loss of generality that
$k\in H^2_{\Cal M}(X_0(\goth n)_L,\Bbb Z(2))$. Let $Y$ be a smooth, projective curve over $\Bbb F_q$ whose function field is $L$ and let $\goth X$ be a regular flat projective model of $X_0(\goth n)_L$ over $Y$. By passing to a finite extension of  $L$, if it is necessary, we may also assume that $\goth X$ is semi-stable. By Proposition 8.6 there is a positive integer $j$ such that the element $jC_*(k)\in H^2_{\Cal M}(X_0(\goth n)_L,\Bbb Z(2))$ lies in the image of the natural map $H^2_{\Cal M}(\goth X,\Bbb Z(2))\rightarrow H^2_{\Cal M}(X_0(\goth n)_L,\Bbb Z(2))$. By Proposition 6.5 of [23] the group $H^2_{\Cal M}(\goth X,\Bbb Z(2))$ is the extension of a torsion group by a $p$-divisible subgroup. The image of the restriction of the regulator of Notation 5.10 to $H^2_{\Cal M}(X_0(\goth n)_L,\Bbb Z(2))$ lies in $\Cal H(\goth n,F^*_{\infty})$ so its image has a torsion $p$-divisible part. Therefore the image of the restriction of this regulator to $H^2_{\Cal M}(\goth X,\Bbb Z(2))$ is torsion. The claim is now clear.\ $\square$
\enddefinition
\definition{Remarks 5.13} Let $\Bbb T(\goth n)$ denote the algebra with unity generated by the endomorphisms $T_{\goth m}$ of the Jacobian $J_0(\goth n)$, where $\goth m\triangleleft A$ is any proper ideal. The algebra $\Bbb T(\goth n)$ is known to be commutative. By claim $(iii)$ of
Proposition 5.11 and the lemma above the algebra of correspondences
generated by the Hecke correspondences leaves the kernel of the regulator of Notation 5.10 restricted to $H^2_{\Cal M}(X_0(\goth n)_L,K(2))$ invariant and its action on the image of this homomorphism
factors through the Hecke al\-geb\-ra $\Bbb T(\goth n)\otimes\Bbb Q$. Moreover the Hecke operator $T_{\goth m}$ acts on this image via the operator $T_{\goth m}$ given by the formula in Definition 5.7 by claim $(iii)$ of Proposition 5.11.
\enddefinition
\definition{Definition 5.14} Let $\mu_G$ be the unique left-invariant Haar measure on the locally compact abelian topological group
$GL_2(\Bbb A)/Z(\Bbb A)$ such that $\mu_G(GL_2(\Cal O)/Z(\Cal O))$ is equal
to $1$. Since this measure is left-invariant with respect to the discrete
action of the group $GL_2(F)/Z(F)$, it induces a measure $Z(\Bbb A)GL_2(F)
\backslash GL_2(\Bbb A)$ which will be denoted by the same
symbol by abuse of notation. Let $V$, $W$ are vector spaces over $\Bbb Q$,
and let $\phi$, $\psi$ be a $V$-valued and a $W$-valued, locally constant
function on $Z(\Bbb A)GL_2(F)\backslash GL_2(\Bbb A)$, respectively. Also assume
that $\psi$ has compact support, for example $\psi\in\Cal H_0(\goth n,W)$.
Then the integral
$$\mathop{\int}_{Z(\Bbb A)GL_2(F)\backslash GL_2(\Bbb A)}
\!\!\!\!\!\!\!\!\!\!\!\!\!\!\!\!\!\!\!\!
\phi(g)\otimes\psi(g)d\mu_G(g)\in V\otimes_{\Bbb Q}W$$
is well-defined. It will be denoted by
$\langle\phi,\psi\rangle$, and will be called the Pe\-tersson product of
$\phi$ and $\psi$.
\enddefinition
\proclaim{Lemma 5.15} For every $k\in H^2_{\Cal M}(Y_0(\goth n)_L,
K(2))$ there is a $k'\in H^2_{\Cal M}(X_0(\goth n)_L,K(2))$ such that
$\langle\{k\},\psi\rangle=\langle\{k'\},\psi\rangle$
for every $\psi\in\Cal H_0(\goth n,\Bbb Q)$.
\endproclaim
\definition{Proof} Let $\Cal U$ denote the group $H^0(Y_0(\goth n)_L,\Cal O^*)$ and let $\Cal V$ denote the $K$-vector subspace of
$H^2_{\Cal M}(Y_0(\goth n)_L,K(2))$ generated by the
product $L^*\otimes\Cal U\subseteq H^2_{\Cal M}(Y_0(\goth n)_L,
\Bbb Z(2))$. By our assumptions on $F$ and $\infty$ the curve
$X_0(\goth n)$ is geometrically irreducible. Moreover 
the group generated by the linear equivalence class of degree zero
divisors defined over $L$ supported on the complement of $Y_0(\goth n)_L$ in the Jacobian of $X_0(\goth n)_L$ is finite by the main theorem of [8]. Hence there is a $u\in\Cal V$ such that $\{k\}_x=\{u\}_x$ for every
closed point $x$ in the complement of $Y_0(\goth n)_L$ where $\{\cdot\}_x$ denotes the $K$-linear extension of the tame symbol at $x$ (this fact is referred to as Bloch's lemma in [24]). Therefore
$k'=k-u$ lies in $H^2_{\Cal M}(X_0(\goth n)_L,K(2))$ by the
localization sequence. Hence it will be sufficient to prove that
$\langle\{u\},\psi\rangle=0$ for every $u\in\Cal V$ and for every
$\psi\in\Cal H_0(\goth n,\Bbb Q)$. In fact we will show the same
claim for every $\psi\in\Cal H_0(\goth n,\overline{\Bbb Q})$. The
operators $T_{\goth m}$ act semisimply on the finite-dimensional
vector space $\Cal H_0(\goth n,\overline{\Bbb Q})$ therefore the
latter decomposes as the direct sum of Hecke-eigenspaces. Hence
we may assume that $\psi$ above is a Hecke-eigenform by linearity.
By the projection formula for the norm map in Milnor $K$-theory we
have $T_{\goth m}(u_1\otimes u_2)=u_1\otimes T_{\goth m}(u_2)$
for every $u_1\in L^*$, $u_2\in\Cal U$ and Hecke correspondence
$T_{\goth m}$. Let $\goth q\!\!\!\not|\goth n$ be a non-zero prime ideal which has a generator $\pi\in A$ such that $\pi\equiv1\mod\goth n$. Then the Hecke correspondence $T_{\goth q}$
maps the cusps, the geometric points in the complement of $Y_0(\goth n)$, into themselves with multiplicity $1+q^{\deg(\goth q)}$ according to the proof of Proposition 3.1 of [8] on page 365. (Strictly speaking this claim is proved for the cusps of the Drinfeld modular curve $X(\goth n)$ there but the former claim immediately follows from the latter.) Hence we have:
$$(1+q^{\deg(\goth q)})\langle\{u\},\psi\rangle=
\langle T_{\goth q}\{u\},\psi\rangle=
\langle\{u\},T_{\goth q}(\psi)\rangle=
\psi^*(\goth q)\langle\{u\},\psi\rangle$$
using the self-adjointness of the operator $T_{\goth q}$
with respect to the Petersson product where
$\psi^*(\goth q)\in\overline{\Bbb Q}$ is the $\goth q$-th Hecke
eigenvalue of $\psi$. By the Ramanujan-Petersson conjecture (proved
in [4] first in this case) the latter is not equal to
$1+q^{\deg(\goth q)}$ when $\deg(\goth q)$ is sufficiently large
hence $\langle\{u\},\psi\rangle$ must be zero.\ $\square$
\enddefinition

\heading 6. The Rankin-Selberg method
\endheading

\definition{Notation 6.1} Let $\goth m\triangleleft A$ be a proper ideal. 
Recall that a Dirichlet character of conductor $\goth m$ is a continuous homomorphism $\chi:\Bbb A_f^*\rightarrow\Bbb C^*$ which is trivial on $F^*\Cal O_{\goth m}$ where $\Cal O_{\goth m}=
\{u\in\Cal O_f^*|u\equiv1\!\!\!\mod\goth m\}$. Then there is a unique homomorphism from $(\Cal O_f/\goth m)^*$ into $\Bbb C^*$, which 
will be denoted by $\chi$ as well by abuse of notation, such that the latter is trivial on the class of constants $\Bbb F_q^*\subset\Cal O_f^*$ and we have $\chi(z)=\chi(\overline z)$ for every $z\in\Cal O^*_f$ where $\overline z$ denote the class of $z$ in the quotient group
$(\Cal O_f/\goth m)^*$. Moreover we let $\chi$ denote also the unique extension of these two homomorphisms onto $\Bbb A_f$ and
$\Cal O_f/\goth m$ which is zero on the complement of $\Bbb A_f^*$ and $(\Cal O_f/\goth m)^*$, respectively.  We are going to assume that the homomorphism $\chi$ is non-trivial. In this case 
$\sum_{\alpha\in\Cal O_f/\goth m}\chi(\alpha)$ is zero. Let $\chi_1$,
$\chi_2\in\Bbb C[\Cal V_{\goth m}]_0$ denote the functions given by the rules:
$$\chi_1(\alpha,\beta)=\chi(\alpha)\quad(\forall(\alpha,\beta)\in\Cal V_{\goth m})
\text{, and}$$
$$\chi_2(\alpha,\beta)=\chi(\beta)\quad\text{when $\alpha=0$,}$$
and $\chi_2(\alpha,\beta)=0$, otherwise, respectively. Let $E^G_{\goth m}(\chi,g,x,y)$ denote the function $\chi(\det(g_f))^{-1}
E^G_{\goth m}(\chi_1,\chi_2,g,x,y)$ for every finite quotient $G$ of $F^*\backslash
\Bbb A^*/\Cal O_f^*$. 
\enddefinition
\proclaim{Lemma 6.2} The function $E^G_{\goth m}(\chi,g,x,y)$ is
left-invariant with respect to $GL_2(F)$ and right-invariant with respect to $\Bbb K(\goth m\infty)\Gamma_{\infty}Z(\Bbb A)$.
\endproclaim
\definition{Proof} By claim $(i)$ of Proposition 3.3 we only need to show that $E^G_{\goth m}(\chi,g,x,y)$ is right-invariant with respect to
$Z(\Bbb A_f)$.  But $Z(\Bbb A_f)=Z(F)Z(\Cal O_f)$ hence
we only have to show that $E^G_{\goth m}(\chi,g,x,y)$ is right-invariant with respect to $Z(\Cal O_f)$. In order to do so we will introduce some convenient notation which we will also use later on without further notice. By our usual abuse of notation for $i=1$, $2$ let $\chi_i:
\Cal O_f^2\rightarrow\Bbb C$ denote the function such that $\chi_i(f)=
\chi_i(\overline f)$ for every $f\in\Cal O_f^2$ where $\overline f$ denotes the class of $f$ in the quotient group $(\Cal O_f/\goth m)^2$. For every $g\in GL_2(\Bbb A)$ let
$$W(g)=\{0\neq f\in F^2|fg_f\in\Cal O_f^2\}\text{,\ }
V(g)=\{f\in W(g)|fg_{\infty}\in F^2_<\}$$
and $U(g)=W(g)-V(g)$. For every $g\in GL_2(\Bbb A)$ and
$z\in Z(\Cal O_f)=\Cal O_f^*$ we have $U(gz)=U(g)$ and $V(gz)=V(g)$. Moreover we have $\chi_i(fz)=\chi_i(f)\chi(z)$ for every $f\in\Cal O_f^2$ and $i=1$, $2$. Therefore
$$\split E^G_{\goth m}(\chi,gz,x,y)=&\chi(z)^{-2}\chi(\det(g_f))^{-1}
{\det(z^{-1})^G\det(g_f^{-1})^G\over(xy)^{\deg(\det(z))+\deg(\det(g))}}\\
&\cdot\!\!\!\!\!\!\!\!
\sum\Sb(a,b)\in U(g)\\(c,d)\in
V(g)\endSb\!\!\!\!\bigg(\chi_1((a,b)g_fz)\chi_2((c,d)g_fz)\\
&\cdot\det\left(\matrix a&\!\!b\\c&\!\!d\endmatrix\right)_{\infty}^G\!\!\!\!\!
x^{2\infty((a,b)g_{\infty})}y^{2\infty((c,d)g_{\infty})}\bigg)\\
&=E^G_{\goth m}(\chi,g,x,y),\endsplit$$
because $\deg(\det(z))=0$ and $\det(z)^G=1$
by definition.\ $\square$
\enddefinition
\definition{Definition 6.3} Let $\chi_0:\Cal O_f\rightarrow\Bbb C$ denote the function such that $\chi_0(u)=\chi(\overline u)$ for every $u\in\Cal O_f$ where $\overline u$ denotes again the class of $u$ in the quotient group $\Cal O_f/\goth m$.  Let $G$ be a finite quotient group of $F^*\backslash\Bbb A^*/\Cal O_f^*$. Note that for every
non-zero $\goth q\triangleleft A$ the value  $y^G\in G$ depends only on $\goth q$ for every $y\in\Bbb A_f^*$ where the divisor of $y$ is
$\goth q$. Let $\goth q^G$ denote this common value. Similarly note that for every non-zero ideal $\goth q\triangleleft A$ relatively prime to $\goth m$ the value $\chi_0(a)$ depends only on $\goth q$ for every $a\in A$ which
generates the ideal $\goth q$. We let $\chi(\goth q)$ denote this common value. For every $G$ as above let $L^G_{\goth m}(\chi,x)$ be the infinite series:
$$\sum_{(\goth q,\goth m)=1}\chi(\goth q)(\goth q^G)^{-1}
x^{\deg(\goth q)}\in\Bbb C[G][[x]].$$
Note that for every complex number $s$ the $\Bbb C[G]$-valued series $L^G_{\goth m}(\chi,q^{-s})$ is absolutely convergent when Re$(s)>1$.
\enddefinition
For every $z\in\Bbb A_f^*$ let $L^G_{\goth m}(\chi,z,s)$ denote the $\Bbb C[G]$-valued series:
$$L^G_{\goth m}(\chi,z,s)=
(z^{-1})^G|z|^s\sum\Sb u\in F^*\\uz\in\Cal O_f\endSb
\chi_0(uz)u_{\infty}^G|u|_{\infty}^{-s}$$
if the latter is absolutely convergent. 
\proclaim{Lemma 6.4} For every $z\in\Bbb A_f^*$ we have
$L^G_{\goth m}(\chi,z,s)=\chi(z)(q-1)L^G_{\goth m}(\chi,q^{-s})$.
\endproclaim
\definition{Proof} First we are going to show that the function
$\overline{\chi(z)}L^G_{\goth m}(\chi,s)$ is invariant with respect to $\Cal O_f^*$. For every $v\in\Cal O_f^*$ we have
$$\split\overline{\chi(zv)}L^G_{\goth m}(\chi,\eta zv,s)=&
\overline{\chi(v)}(v^{-1})^G|v|^s
\overline{\chi(z)}(z^{-1})^G|z|^s\cdot\!\!\!\!\!
\sum\Sb u\in F^*\\uzv\in\Cal O_f\endSb
\!\!\!\chi_0(uz)\chi_0(v)
u_{\infty}^G|u|_{\infty}^{-s}\\
=&L^G_{\goth m}(\chi,z,s),\endsplit$$
because $\chi(v)=\chi_0(v)$, $v^G=1$ and $|v|=1$. Now we may assume that $z\in F^*$ because $\Bbb A_F^*=F^*\Cal O_f^*$. In this case
$$\split\overline{\chi(z)}L^G_{\goth m}(\chi,z,s)=L^G_{\goth m}(\chi,z,s)=&
\sum_{0\neq u\in A}\chi_0(u)(u_f^G)^{-1}
|u|^{-s}_{\infty}\\=&(q-1)L^G_{\goth m}(\chi,q^{-s})\endsplit$$
because we have $(u^G_f)^{-1}=u_{\infty}^G$ for every $u\in F^*$, and because $\chi(z)=1$ and $|z_f|=|z|_{\infty}^{-1}$ as the degree of every principal divisor is zero.\ $\square$
\enddefinition
\definition{Definition 6.5} Let $B$ denote the group scheme of invertible upper triangular two by two matrices. For every finite quotient $G$ as above and
$g\in B(\Bbb A)$ let $K^G_{\goth m}(\chi,g,s)$ denote the $\Bbb C[G]$-valued function:
$$K^G_{\goth m}(\chi,g,s)=\overline{\chi((xz)_f)}
((xz)_f^{-1})^G|xz^2|^s\!\!\!\!\!\!\!
\sum_{(v,w)\in U(g)}\!\!\!\!\!\!
\chi_0((vxz)_f)v_{\infty}^G\;|(vxz)_{\infty}|^{-2s}$$
where $g=\left(\smallmatrix xz&yz\\0&z\endsmallmatrix\right)\in B(\Bbb A)$ and $s$ is a complex number when this infinite sum is absolutely convergent. The latter holds if $s$ with Re$(s)>1$ because the series above is majorated by the series $E(g,s)$. Finally for every pair of complex numbers $s$, $t$ with Re$(s)>1$, Re$(t)>1$ we let $H^G_{\goth m}(\chi,g,s,t)$ denote the 
$\Bbb C[G]$-valued function on $GL_2 (\Bbb A)$ given by the formula:
$$H^G_{\goth m}(\chi,g,s,t)=L^G_{\goth m}(\chi,q^{-2t})|x|^t
K^G_{\goth m}(\chi,\left(\matrix xz&yz\\0&z\endmatrix\right),s)$$
if
$$g=\left(\matrix xz&yz\\0&z\endmatrix\right)k\text{\ where\ }
k\in\Bbb K(\goth m\infty)\Gamma_{\infty},$$
and $H^G_{\goth m}(\chi,g,s,t)=0$, otherwise.
\enddefinition
\proclaim{Lemma 6.6} The following holds:
\roster
\item"$(i)$" the function $K^G_{\goth m}(\chi,g,s)$ is
left-invariant with respect to $B(F)$ and right-invariant with respect to
$(B(\Cal O)\cap\Bbb K(\goth m))Z(\Bbb A)$,
\item"$(iii)$" the function $H^G_{\goth m}(\chi,g,s,t)$ is well-defined and it
is left-invariant with respect to $B(F)$ and right-invariant with respect to $\Bbb K(\goth m\infty)\Gamma_{\infty}Z(\Bbb A)$.
\endroster
\endproclaim
\definition{Proof} The proof of the first claim is the same as the proofs of claim $(i)$ of Proposition 3.3 and Lemma 6.4. In order to prove that
$H^G_{\goth m}(\chi,g,s,t)$ is well-defined we need to show that $|x|^t=|a|^t$ and
$$K^G_{\goth m}(\chi,\left(\matrix xz&yz\\0&z\endmatrix\right),s)=
K^G_{\goth m}(\chi,\left(\matrix ac&bc\\0&c\endmatrix\right),s)$$
where
$$\left(\matrix xz&yz\\0&z\endmatrix\right),
\left(\matrix ac&bc\\0&c\endmatrix\right)\in B(\Bbb A)\text{\ and\ }
\left(\matrix xz&yz\\0&z\endmatrix\right)\cdot
\left(\matrix ac&bc\\0&c\endmatrix\right)^{-1}
\in\Bbb K(\goth m\infty)\Gamma_{\infty}.$$
The latter is an immediate consequence of claim $(i)$. Similarly the invariance properties of $H^G_{\goth m}(\chi,g,s,t)$ claimed above are 
obvious from claim $(i)$ and its definition.\ $\square$
\enddefinition
\proclaim{Proposition 6.7} For every $g\in GL_2(\Bbb A)$ the sum on the right hand side below
is absolutely convergent and we have:
$$E^G_{\goth m}(\chi,g,q^{-s},q^{-t})=(q-1)\!\!\!\!\!\!\!\!\!
\sum_{\rho\in B(F)\backslash GL_2(F)}\!\!\!\!\!\!\!\!
H^G_{\goth m}(\chi,\rho g,s,t),$$
when $\text{\rm Re}(s)>1$ and $\text{\rm Re}(t)>1$.
\endproclaim
\definition{Proof} For every $\rho\in GL_2(F)$ the value of
$H^G_{\goth m}(\chi,\rho g,s,t)$ depends only on the left $B(F)$-coset of
$\rho$ because $H^G_{\goth m}(\chi,g,s,t)$ is left-invariant with respect to
$B(F)$. Hence the infinite sum on the right hand side above is well-defined. By grouping the terms of the absolutely convergent series on the left hand side we get:
$$\split E^G_{\goth m}(\chi,g,q^{-s},q^{-t})=&
\sum_{\rho\in B(F)\backslash GL_2(F)}\!\!\!\!\!\!\!\!\!\!
\chi(\det((\rho g)_f))^{-1}
\det((\rho g)_f^{-1})^G|\det(\rho g)|^{s+t}\\
&\cdot\Big(\sum\Sb(v,w)\in U(\rho g)
\\u\in F^*,(0,u)\in V(\rho g)\endSb
\!\!\!\!\!\!\!\!\!\!\!\!
\chi_1((v,w)\rho g_f)\chi_2((0,u)\rho g_f)(vu)_{\infty}^G\\
&\cdot
\|(v,w)\rho g_{\infty}\|^{-2s}\cdot\|(0,u)\rho g_{\infty}\|^{-2t}\Big).
\endsplit$$
By the Iwasawa decomposition we may write $g$ as
$$g=pk\text{, where\ }p=\left(\matrix a&b\\0&c\endmatrix\right)
\in B(\Bbb A)\text{\ and\ }
k=\left(\matrix k_{11}&k_{12}\\k_{21}&k_{22}\endmatrix\right)
\in GL_2(\Cal O).$$
Because $k_{\infty}$ is an isometry we only have to show that
$$\split H^G_{\goth m}(\chi,g,s,t)=&
\Big(\overline{\chi(c_f)}(c_f^{-1})^G|ac|^t\!\!\!\!\!\!
\sum\Sb u\in F^*\\(0,u)\in V(g)\endSb\!\!\!\!\!\!\!
\chi_2((0,u)g_f)u^G_{\infty}\|(0,u)p_{\infty}\|^{-2t}\Big)\\&\cdot
\Big(\overline{\chi(a_f)}(a_f^{-1})^G|ac|^s\!\!\!\!\!\!\!
\sum_{(v,w)\in U(g)}\!\!\!\!\!\!\!
\chi_1((v,w)g_f)v_{\infty}^G\|(v,w)p_{\infty}\|)^{-2s}\Big)\endsplit$$
by the above. The first infinite sum is zero unless there is a
$d\in\Bbb A_f$ such that $d(k_{21},k_{22})_f\in\Cal O_f^2$ and the latter is congruent to $(0,\alpha)$ modulo $\goth m\Cal O_f$ for some $\alpha\in(\Cal O_f/\goth m)^*$. The latter is possible exactly when $k_f$ is in $\Bbb K_0(\goth m)$. We may even assume that $k_f$ is in $\Bbb K(\goth m)$ by changing $p$, if necessary. By the definition of the set $V(g)$ we  also need that
$|(k_{22})_{\infty}|>|(k_{21})_{\infty}|$ for the first sum to be non-zero. Since $(k_{22})_{\infty}\in\Cal O_{\infty}$ we have
$\infty((k_{21})_{\infty})>0$ so $k_{\infty}\in\Gamma_{\infty}$. In this case we have
$(0,u)g_{\infty}\in F^2_<$ for every $u\in F^*$ automatically so
$$\{u\in F^*|(0,u)\in V(g)\}=
\{u\in F^*|uc\in\Cal O_f\}.$$
Hence the first term of the product above is
$|a/c|^t\overline{\chi(c_f)}L^G(\overline\chi,c_f, 2t)$. By Lemma 6.4
we know that  the latter is $|a/c|^t(q-1)L(\overline\chi,q^{-2t})$. On the other hand the second term is visibly $K^G_{\goth m}(\chi,g,s)$ because $U(g)=U(p)$ and $\chi_1(fk_f)=\chi_1(f)$ for every $f\in\Cal O_f^2$ since $k_f\in\Bbb K(\goth m)$.\ $\square$
\enddefinition
\definition{Definition 6.8} Let $\mu_B$ be the unique left-invariant Haar measure on the locally compact abelian topological group $Z(\Bbb A)\backslash B(\Bbb A)$ such that $\mu_B(Z(\Cal O)\backslash B(\Cal O))$ is equal to $1$.  Since this measure is left-invariant with respect to the discrete action of the group $Z(F)\backslash B(F)$, it induces a measure on $Z(\Bbb A)B(F)\backslash B(\Bbb A)$, which will be denoted by the same symbol by abuse of notation. The measure $\mu_B$ has a simple description. Let $\mu$ and $\mu^*$ be the unique Haar measure on the locally compact abelian
topological group $\Bbb A$ and $\Bbb A^*$, respectively, such that
$\mu(\Cal O)$ and $\mu^*(\Cal O^*)$ are both equal to $1$.  Since the measures $\mu$ and $\mu^*$ are
left-invariant with respect to the discrete subgroups $F\subset\Bbb A$,
and $F^*\subset\Bbb A^*$, respectively, by definition, they induce a measure on $F\backslash\Bbb A$ and $F^*\backslash\Bbb A^*$, respectively,  which will be denoted by the same letter by abuse of notation. Then we have
$$\mathop{\int}_{Z(\Bbb A)B(F)\backslash B(\Bbb A)}
\!\!\!\!\!\!\!\!\!\!\!\!\!\!\!
f\left(\matrix x&y\\0&1\endmatrix\right)\mu_B(\left(\matrix x&y\\0&1\endmatrix\right))
=
\!\!\!\!\!\mathop{\int}_{F^*\backslash\Bbb A^*}\!\!\!\!\!\mu^*(x)
\!\!\!\mathop{\int}_{F\backslash\Bbb A}\!\!\!
f\left(\matrix x&y\\0&1\endmatrix\right){\mu(y)\over|x|}$$
for every Lebesgue-measurable function $f:Z(\Bbb A)B(F)\backslash B(\Bbb A)\rightarrow\Bbb C$.
\enddefinition
\proclaim{Lemma 6.9} For every
$\psi\in\Cal A_0(\goth m\infty,\Bbb C)$ the integrands of
the two integrals below are absolutely Lebesgue-integrable and
$$\mathop{\int}_{Z(\Bbb A)GL_2(F)\backslash GL_2(\Bbb A)}
\!\!\!\!\!\!\!\!\!\!\!\!\!\!\!\!\!\!\!\!
E^G_{\goth m}(\chi,g,q^{-s},q^{-t})
\overline{\psi(g)}d\mu_G(g)
=\mu(\goth m)\!\!\!\!\!\!\!\!\!\!\!\!\!\!\!
\mathop{\int}_{Z(\Bbb A)B(F)\backslash B(\Bbb A)}
\!\!\!\!\!\!\!\!\!\!\!\!\!\!\!
H^G_{\goth m}(\chi,b,s,t)\overline{\psi(b)}d\mu_B(b)$$
where $\mu(\goth m)=(q-1)
\mu_G(Z(\Cal O)\backslash\Bbb K(\goth m\infty)
\Gamma_{\infty}Z(\Cal O))$ when
$\text{\rm Re}(s)>1$ and $\text{\rm Re}(t)>1$.
\endproclaim
\definition{Proof} We may talk about the Lebesgue-integrability of the integrands above because they are $\Bbb C[G]$-valued functions. By Theorem 2.2.1 in [11], pages 255-256, we know that any cuspidal automorphic form which is invariant with respect to
$Z(\Bbb A)$ has compact support as a function on $Z(\Bbb A)GL_2(F)
\backslash GL_2(\Bbb A)$. Hence the integral on the left in the
equation above is absolutely convergent and we may interchange the
integration and the summation in Proposition 6.7 to get that
$$\mathop{\int}_{Z(\Bbb A)GL_2(F)\backslash GL_2(\Bbb A)}
\!\!\!\!\!\!\!\!\!\!\!\!\!\!\!\!\!\!\!\!
E^G_{\goth m}(\chi,g,q^{-s},q^{-t})
\overline{\psi(g)}d\mu_G(g)
=(q-1)\!\!\!\!\!\!\!\!\!\!\!\!\!\!\!
\mathop{\int}_{Z(\Bbb A)B(F)\backslash GL_2(\Bbb A)}
\!\!\!\!\!\!\!\!\!\!\!\!\!\!\!\!\!
H^G_{\goth m}(\chi,b,s,t)\overline{\psi(b)}d\mu_G(b)$$
where the measure on $Z(\Bbb A)B(F)\backslash GL_2(\Bbb A)$ induced by $\mu_G$ will be denoted by the same symbol by the usual abuse of notation. The map:
$$\pi:Z(\Bbb A)B(F)\backslash B(\Bbb A)\times
Z(\Cal O)\backslash\Bbb K(\goth m\infty)\Gamma_{\infty}
Z(\Cal O)
\rightarrow Z(\Bbb A)B(F)\backslash GL_2(\Bbb A)$$
given by the rule $(b,k)\mapsto bk$ is continuous, hence for every
Borel-measurable set $\goth B\subseteq Z(\Bbb A)B(F)\backslash GL_2(\Bbb A)$ the pre-image $\pi^{-1}(\goth B)$ is also Borel-measurable. Let $\mu_B\times\mu_G$ denote the direct product of the measures $\mu_B$ and $\mu_G$ on the direct product $Z(\Bbb A)B(F)\backslash B(\Bbb A)\times Z(\Cal O)\backslash
\Bbb K(\goth m\infty)\Gamma_{\infty}Z(\Cal O)$. Then we have
$\mu_B\times\mu_G(\pi^{-1}(\goth B))=\mu_G(\goth B)$ for every $\goth B$ above. Moreover the map $\pi$ maps surjectively onto the support of  $H^G_{\goth m}(\chi,b,s,t)$ as a function on $GL_2(\Bbb A)$
as we saw in Definition 6.4 so the integral above is equal to:
$$(q-1)\!\!\!\!\!\!\!\!\!\!
\mathop{\int}_{Z(\Bbb A)B(F)\backslash B(\Bbb A)}
\!\!\!\!\!\!\!\!\!\!\!\!\!\!d\mu_B(b)
\mathop{\int}_{Z(\Cal O)\backslash\Bbb K(\goth m\infty)
\Gamma_{\infty}Z(\Cal O)}\!\!\!\!\!\!\!\!\!\!\!\!
H^G_{\goth m}(\chi,bk,s,t)\overline{\psi(bk)}d\mu_G(k)$$
by Fubini's theorem. By definition 
the integrand of the interior integral is constant on the
domain of integration. The claim is now obvious.\ $\square$
\enddefinition
\definition{Definition 6.10} Let $\tau:F\backslash\Bbb A\rightarrow\Bbb C^*$ be a non-trivial continuous character and let $\goth d$ be an idele such that $\Cal D=\goth d\Cal O$, where $\Cal D$ is the $\Cal O$-module defined as
$$\Cal D=\{x\in\Bbb A|\tau(x\Cal O)=1\}.$$
It is well-known the linear equivalence class of the divisor of $\goth d$ is the anti-canonical class. Moreover for every $\eta\in F^*$ the map
$x\mapsto\tau(\eta x)$ is another non-trivial continuous homomorphism. Therefore by choosing an appropriate character $\tau$, we may assume that $\goth d$ is any idele of degree two, as every such divisor is linearly equivalent to the anti-canonical class. In particular we may assume that $\goth d=\pi^2$
where $\pi\in F_{\infty}$ is a uniformizer. For every $\goth r\triangleleft A$ non-zero ideal let
$S(\goth m,\goth r)$ denote the set:
$$S(\goth m,\goth r)=
\{0\neq\goth q\triangleleft A|(\goth m,\goth q)=1,\goth q|\goth r\}.$$
Moreover for every $G$ as in Definition 6.3 let $\sigma^G_{\goth m}(\chi,\goth r,x)\in\Bbb C[G][x]$ denote the polynomial given by the formula:
$$\sigma^G_{\goth m}(\chi,\goth r,x)=\!\!\!\!
\sum_{\goth q\in S(\goth m,\goth r)}\chi(\goth q)
(\goth q^G)^{-1}x^{\deg(\goth q)}.$$
\enddefinition
\proclaim{Proposition 6.11} For each complex $s$ with $\text{\rm Re}(s)>1$
we have:
$$\mathop{\int}_{F\backslash\Bbb A}
K^G_{\goth m}(\chi,\left(\matrix x\goth d&y\\0&1\endmatrix\right),s)
\tau(-y)d\mu(y)=(q-1)
|x\goth d|^{1-s}
\sigma^G_{\goth m}(\chi,x_f,q^{1-2s}),$$
if the divisor of $x$ is effective, and it is zero, otherwise.
\endproclaim
\definition{Proof} The integral above is well-defined because the integrand is $F$-invariant by claim $(i)$ of Lemma 6.6. Note that we have $u\neq 0$ for every $0\neq(u,v)\in F^2$ such that
$\chi_0((ux\goth d)_f)\neq0$. Therefore by grouping the terms in the
infinite sum of Definition 6.5 we get the following identity:
$$K^G_{\goth m}(\chi,\!\left(\matrix x\goth d&\!\!\!y\\0&\!\!\!1\endmatrix
\right)\!,s)=
\overline{\chi((x\goth d)_f)}((x\goth d)_f^{-1})^G|x\goth d|^s
\sum_{v\in F}\!\!\!\!\!\!\!\!\!\!\!\!\!\!\!\!\!\!
\sum\Sb u\in F^*\\(u,0)\in U(
\left(\smallmatrix x\goth d&(y+v)\\0&1\endsmallmatrix\right))\endSb
\!\!\!\!\!\!\!\!\!\!\!\!\!\!\!\!\!\!
\chi_0((ux\goth d)_f)u_{\infty}^G|(ux\goth d)_{\infty}|^{-2s}.$$
Hence
$$\split\mathop{\int}_{F\backslash\Bbb A}
K^G_{\goth m}(\chi,\left(\matrix x\goth d&y\\0&1\endmatrix\right),s)
\tau(-y)d\mu(y)=&\\
\overline{\chi((x\goth d)_f)}
((x\goth d)_f^{-1})^G\cdot|x\goth d|^s\!\!\!
&\!\!\!\!\sum\Sb u\in F^*\\u(x\goth d)_f\in\Cal O_f\endSb\!\!\!\!\!\!
\chi_0((ux\goth d)_f)u_{\infty}^G|ux\goth d|_{\infty}^{-2s}
\!\!\!\!\!\!\!\!\!\!\!
\mathop{\int}\Sb uy_f\in\Cal O_f
\\|y_{\infty}|\leq|x\goth d|_{\infty}\endSb
\!\!\!\!\!\!\!\!\!\!\!\tau(-y)d\mu(y)\endsplit$$
by interchanging summation and integration. For every $u\in F^*$ the domain of integration of the integral above is a
direct product of the sets $u_f^{-1}\Cal O_f\subset\Bbb A_f$ and
$x\goth d\Cal O_{\infty}\subset F_{\infty}$. The integral itself is
non-zero if and only if the product set above lies in the
kernel of $\tau$. The latter is equivalent to the conditions
$(u\goth d)_f^{-1}=u_f^{-1}\in\Cal O_f$ and $\infty(x)\geq0$. In this case
the integral is equal to:
$$\mu(u_f^{-1}\Cal O_f\times x\goth d\Cal O_{\infty})=
|u|^{-1}\mu(\Cal O_f\times ux\goth d\Cal O_{\infty})
=|ux\goth d|_{\infty}.$$
Let $T(\goth m,x)$ denote the set:
$$T(\goth m,x)=\{u\in F^*|(ux)_f\in\Cal O_f,
u_f^{-1}\in\Cal O_f\}.$$
By the above the left hand side of the equation in the claim above is equal to:
$$\overline{\chi(x_f)}
|x\goth d|^{1-s}\!\!\!\sum_{u\in T(\goth m,x)}\!\!\!
\chi_0((ux)_f)((ux)_f^{-1})^G|(ux)_f|^{2s-1}$$
when $\infty(x)\geq0$, and it is zero, otherwise. The set $T(\goth m,x)$ is empty when $x_f$ is not an element of $\Cal O_f$. Therefore the expression above is zero unless the divisor of $x$ is effective. Note that in the latter case for every $u\in T(\goth m,x)$ the number $\chi_0((ux)_f)$ is zero unless the divisor of $(ux)_f$ is an element of $S(\goth m,x_f)$. On the other hand every element of $S(\goth m,x_f)$ is the divisor of an idele of the form $(ux)_f$ for some $u\in T(\goth m,x)$ and $u$ is unique up to factor in $\Bbb F_q^*$. Note that the sum above is invariant in the variable $x$ with respect to the action of $\Cal O^*_f$. Hence we may assume that $x_f=\eta_f$ for some $\eta\in F^*$. In this case we have $\chi(x_f)=1$
and $\chi(\goth q)=\chi_0((ux)_f)$ for every $u\in T(\goth m,x)$ where $\goth q$ is the divisor of $(ux)_f$. The claim is now clear.\ $\square$
\enddefinition
\definition{Notation 6.12} Recall that we call two divisors $\goth r$ and $\goth s$ on $X$ relatively prime if their support is disjoint. For every $\psi\in\Cal A_0(\goth m,\Bbb C)$ let $\psi^*:\text{\rm Div}(X)\rightarrow\Bbb C$ denote the Fourier coefficients of $\psi$ whose existence was established in Proposition 1 of Chapter III in [25], page 21, proved on pages 19-20 of [25]. Recall that a function $f:\text{Div}(X)\rightarrow R$ is called multiplicative, where $R$ is a commutative ring with unity, if it is zero on non-effective divisors, $f(1)=1$ and for every pair of relatively prime divisors $\goth r$ and $\goth s$ we have $f(\goth r\goth s)=f(\goth r)f(\goth s)$. (Similarly an $R$-valued function on the set of non-zero ideals of $A$ is called multiplicative if it satisfies the last two properties of the previous definition.)
Let us recall the situation considered in the introduction. Let $E$ be an elliptic curve defined over $F$ which has split multiplicative reduction at $\infty$. By assumption the conductor of $E$ is
of the form $\goth n\infty$ where $\goth n$ is an effective divisor which is supported in the complement of $\infty$ in $X$. Let $\psi^*_E$ denote the unique multiplicative function into the
multiplicative semigroup of $\Bbb Q$ such that $\psi^*_E(x^n)$ is the
same as in 1.6 for each natural number $n$ and each closed point $x$ on $X$. A cuspidal harmonic form $\phi_E\in\Cal H_0(\goth n,\Bbb Q)$
is called a normalized Hecke eigenform attached to $E$ is 
if its Fourier coefficient $\phi_E^*(\goth q)$ is equal to $|\goth q|
\psi^*_E(\goth q)$ for every effective divisor $\goth q$.
\enddefinition
The following proposition is an easy consequence of the Langlands
correspondence:
\proclaim{Proposition 6.13} There is a unique normalized Hecke eigenform
attached to $E$.
\endproclaim
\definition{Proof} The only not entirely obvious fact is that the
normalized Hecke eigenform has values in $\Bbb Q$, see for example
the proof of Proposition 3.3 in [21].\ $\square$
\enddefinition
\proclaim{Theorem 6.14} Assume that $\goth n$ divides $\goth m$.
Then for every $\text{\rm Re}(s)>1$ and $\text{\rm Re}(t)>1$ we have:
$$\split\mathop{\int}_{Z(\Bbb A)GL_2(F)\backslash GL_2(\Bbb A)}
\!\!\!\!\!\!\!\!\!\!\!\!\!\!\!\!\!\!\!
&E^G_{\goth m}(\chi,g,q^{-s},q^{-t})
\phi_E(g)d\mu_G(g)=\\&
(q-1)\mu(\goth m)L^G_{\goth m}(\chi,q^{-2t})
{|\goth d|^{t-s}\over1-q^{s-t-1}}
\sum_{0\neq\goth r\triangleleft A}
|\goth r|^{1+t-s}
\sigma^G_{\goth m}(\chi,\goth r,q^{1-2s})\psi_E^*(\goth r).\endsplit$$
\endproclaim
\definition{Proof} By Lemma 6.9 and the description of the measure
$\mu_B$ as a double integral at the end of Definition 6.8 we know that the integral on the left hand side of the equation above is equal to:
$$\mu(\goth m)L^G_{\goth m}(\chi,q^{-2t})
\!\!\!\!\!\mathop{\int}_{F^*\backslash\Bbb A^*}\!\!\!\!\!\mu^*(x)
\!\!\!\mathop{\int}_{F\backslash\Bbb A}\!\!\!
|x|^tK^G_{\goth m}(\chi,\left(\matrix x&y\\0&1\endmatrix\right),s)
\phi_E(\left(\matrix x&y\\0&1\endmatrix\right))
d{\mu(y)\over|x|}.\eqno6.14.1$$
Using the Fourier expansion of $\phi_E=\overline{\phi_E}$ we get from Proposition 6.11 that
$$\split\int_{F\backslash\Bbb A}
K^G_{\goth m}(\chi,\left(\matrix x\goth d&y\\0&1\endmatrix\right),s)
&\phi_E
(\left(\matrix x\goth d&y\\0&1\endmatrix\right))d\mu(y)
=\\
&(q-1)\sum_{\eta\in F^*}
|\eta x\goth d|^{1-s}
\sigma^G_{\goth m}(\chi,(\eta x)_f,q^{1-2s})
\phi_E^*(\eta x).\endsplit$$
By plugging the equation above into the double integral in 6.14.1 
we get that the latter is equal to:
$$(q-1)\!\!\mathop{\int}_{\Bbb A^*}
|x\goth d|^{t-s}
\sigma^G_{\goth m}(\chi,x_f,q^{1-2s})\phi_E^*(x)d\mu^*(x)$$
if we also interchange the summation in the index $\eta$ and the integration. The integrand above is constant on the cosets of the subgroup $\Cal O^*\subset\Bbb A^*$ hence the integral is equal to the infinite sum:
$$\split|\goth d|^{t-s}\sum\Sb0\neq\goth r\triangleleft A\\k\in\Bbb N\endSb
|\goth r\infty^k|^{t-s}
\sigma^G_{\goth m}(\chi,\goth r,q^{1-2s})&
\phi_E^*(\goth r\infty^k)=\\&
{|\goth d|^{t-s}\over1-q^{s-t-1}}
\sum_{0\neq\goth r\triangleleft A}
|\goth r|^{1+t-s}
\sigma^G_{\goth m}(\chi,\goth r,q^{1-2s})\psi_E^*(\goth r),\endsplit$$
where we also used that the function $\psi_E^*$ is multiplicative and
$\psi_E^*(\infty^k)=1$ for every $k\in\Bbb N$.\ $\square$
\enddefinition

\heading 7. An $\infty$-adic analogue of Beilinson's theorem
\endheading

Let $R$ be an arbitrary commutative ring with unity. Let $*:
R[[t]]\times R[[t]]\rightarrow R[[t]]$ denote the map given by the
rule:
$$(\sum_{n\in\Bbb N}a_nt^n)*(\sum_{n\in\Bbb N}b_nt^n)=
\sum_{n\in\Bbb N}a_nb_nt^n.$$
\proclaim{Lemma 7.1} We have:
$$\split{1\over(1-\alpha_1t)(1-\beta_1t)}*&
{1\over(1-\alpha_2t)(1-\beta_2t)}=\\
&{1-\alpha_1\beta_1\alpha_2\beta_2t^2
\over(1-\alpha_1\alpha_2t)(1-\alpha_1\beta_2t)
(1-\beta_1\alpha_2t)(1-\beta_1\beta_2t)}\endsplit$$
for every $\alpha_1$, $\beta_1\in R$ and $\alpha_2$, $\beta_2\in R$.
\endproclaim
\definition{Proof} We may assume that $R=\Bbb Z[x_1,y_1,x_2,y_2]$
and $\alpha_i=x_i$, $\beta_i=y_i$ for $i=1$, $2$ without the loss of
generality. Note that
$${\alpha_i-\beta_i\over(1-\alpha_it)(1-\beta_it)}=
{\alpha_i\over(1-\alpha_it)}-{\beta_i\over(1-\beta_it)}$$
for $i=1$ and $i=2$. Also note that
$${1\over1-\gamma_1t}*{1\over1-\gamma_2t}=
{1\over1-\gamma_1\gamma_2t}$$
for every $\gamma_1$, $\gamma_2\in R$ by definition. Because the map $*$ is $R$-bilinear we have
$$\split(\alpha_1-\beta_1)(\alpha_2-\beta_2)
{1\over(1-\alpha_1t)(1-\beta_1t)}\!\!\!\!\!\!&\,\,\,\,\,\,*
{1\over(1-\alpha_2t)(1-\beta_2t)}\\
\!\!\!\!\!\!=&\left({\alpha_1\over1-\alpha_1t}-{\beta_1\over1-\beta_1t}\right)*
\left({\alpha_2\over1-\alpha_2t}-{\beta_2\over1-\beta_2t}\right)\\
\!\!\!\!\!\!=&
{(\alpha_1-\beta_1)(\alpha_2-\beta_2)(1-\alpha_1\beta_1\alpha_2
\beta_2t^2)
\over(1-\alpha_1\alpha_2t)(1-\alpha_1\beta_2t)
(1-\beta_1\alpha_2t)(1-\beta_1\beta_2t)}\endsplit$$
by the above. Since our assumption above implies that $\alpha_i-\beta_i$ is not a zero divisor in $R$ for $i=1$, $2$ the claim is now
clear.\ $\square$
\enddefinition
\definition{Notation 7.2} Let us consider the situation described in Definition 1.6. The Galois representation $\chi$ corresponds to a Dirichlet character of conductor $\goth m$ described in Definition 6.1 by class field theory if an embedding of $K$ into the field of complex numbers is also provided. We let $\chi$ denote this Dirichlet character, too. Moreover the profinite completion of the group $F^*\backslash\Bbb A^*/\Cal O_f^*$ and $G_{\infty}$ are canonically isomorphic by class field theory. In particular there is a bijective correspondence between the finite quotients of these groups. These two sets are going to be identified in all that follows. For every effective divisor $\goth d$ on
$X$ let $L_{\goth d}(E,x)$ be the $L$-function:
$$L_{\goth d}(E,t)=L(X(\goth d),\rho,t)\in\Bbb C[t]$$
where we continue to use the notation introduced in the proof of Proposition 2.4.
\enddefinition
\proclaim{Proposition 7.3} We have:
$${L_{\goth m\infty}(E,t)\Cal L^G_{\goth m}(E,\chi,xt)
\over L^G_{\goth m}(\chi,qxt^2)}=
\sum\Sb0\neq\goth q\triangleleft A\\
(\goth q,\goth m)=1\endSb
\psi^*_E(\goth q)\sigma^G_{\goth m}(\chi,\goth q,x)t^{\deg(\goth q)}.$$
\endproclaim
\definition{Proof} Note that the $l$-adic representation
$\rho$ is unramified at every prime ideal $\goth q\triangleleft A$ which does not divide $\goth m$ therefore the local factor $L_{\goth q}(E,t)$
of the Hasse-Weil $L$-function of $E$ at $\goth q$ can be written as
$$L_{\goth q}(E,t)={1\over(1-\alpha(\goth q)t^{\deg(\goth q)})
(1-\beta(\goth q)t^{\deg(\goth q)})},$$
where $\alpha(\goth q)$ and $\beta(\goth q)$ are complex numbers such that $\alpha(\goth q)+\beta(\goth q)=\psi^*_E(\goth q)$ and
$\alpha(\goth q)\cdot\beta(\goth q)=q^{\deg(\goth q)}$. On the other hand it is clear from the definition of
$\sigma^G_{\goth m}(\chi,\goth q,x)$ that the latter is a $K[G][x]$-valued multiplicative function on the set of non-zero ideals of $A$. Therefore the power series in both sides of the equation in the claim above are Euler products, that is, the left hand side and the right hand side of the equation above is equal to:
$$\prod\Sb \goth q\in|X|\\\goth q\not\in supp(\goth m\infty)\endSb
\!\!\!\!\!\!\!\!\!A_{\goth q}(x,y)
\quad\quad\text{ and }
\prod\Sb\goth q\in|X|\\\goth q\not\in supp(\goth m\infty)\endSb
\!\!\!\!\!\!\!\!\!B_{\goth q}(x,y),$$ respectively, where
$$\split A_{\goth q}(x,y)=&{1-\chi(\goth q)(\goth q^G)^{-1}
(qxt^2)^{\deg(\goth q)}\over 
(1-\alpha(\goth q)t^{\deg(\goth q)})
(1-\beta(\goth q)t^{\deg(\goth q)})}\\
&\cdot
{1\over(1-\alpha(\goth q)\chi(\goth q)(\goth q^G)^{-1}
(xt)^{\deg(\goth q)})
(1-\beta(\goth q)\chi(\goth q)(\goth q^G)^{-1}(xt)^{\deg(\goth q)})}
\endsplit$$
and
$$B_{\goth q}(x,y)=\sum_{n=0}^{\infty}\psi^*_E(\goth q^n)
\sigma^G_{\goth m}(\chi,\goth q^n,x)t^{\deg(\goth q)n}$$
for every ${\goth q}\in|X|$ such that $\goth q
\not\in supp(\goth m\infty)$ by the above. Clearly it is sufficient to prove that for every $\goth q$  the factors of these Euler products at $\goth q$ are equal. But the latter follows at once from Lemma 7.1 and the fact that
$$\sum_{n=0}^{\infty}
\sigma^G_{\goth m}(\chi,\goth q^n,x)t^{\deg(\goth q)n}=
{1\over(1-t^{\deg(\goth q)})
(1-\chi(\goth q)(\goth q^G)^{-1}(xt)^{\deg(\goth q)})}
\text{.\ $\square$}$$
\enddefinition
\proclaim{Theorem 7.4} We have:
$$\langle E^G_{\goth m}(\chi,g,x,y),\phi_E\rangle=
(q-1)\mu(\goth m)\left({x\over y}\right)^2\!\!
L_{\goth m}(E,{y\over qx})\Cal L^G_{\goth m}(E,\chi,xy).$$
\endproclaim
\definition{Proof} By definition both sides of the equation above are elements of the ring $\Bbb C[G][[x,y]][x^{-1},y^{-1}]$. But in fact the left and the right hand sides are elements of the ring $\Bbb C[G][x,y,x^{-1},y^{-1}]$ by Propositions 3.5 and 2.4, respectively. We also know that after we substitute $q^{-s}$ and $q^{-t}$ into $x$ and $y$, respectively, both sides of the equation above become absolutely
convergent when $\text{\rm Re}(s)>1$ and $\text{\rm Re}(t)>1$. Therefore it will be sufficient to prove that they are equal after these substitutions by the unique continuation of holomorphic functions. Since $\goth d$ is the anticanonical class, its degree is two, so the integral $\langle E^G_{\goth m}(\chi,g,q^{-s},q^{-t}),\phi_E\rangle$ can be rewritten as the infinite sum:
$$\split(q-1)\mu(\goth m)L^G_{\goth m}(\chi,q^{-2t})&
{|\goth d|^{t-s}\over1-q^{s-t-1}}
\sum_{0\neq\goth r\triangleleft A}
|\goth r|^{1+t-s}
\sigma^G_{\goth m}(\chi,\goth r,q^{1-2s})\psi_E^*(\goth r)=\\&
(q-1)\mu(\goth m)q^{2s-2t}
L_{\goth m}(E,q^{1+s-t})
\Cal L^G_{\goth m}(E,\chi,q^{-s-t})\endsplit$$
by Theorem 6.14 and Proposition 7.3. The claim is now
clear. \ $\square$
\enddefinition
\definition{Notation 7.5} It is clear from Theorem 5.4 that
the irreducible components of the curve $Y(\goth m)_{F_{\infty}}$
are in a bijective correspondence with the set:
$$GL_2(F)\backslash GL_2(\Bbb A_f)/\Bbb K_f(\goth m)$$
of double cosets. In fact for a double coset represented by an element
$g\in GL_2(\Bbb A_f)$ the corresponding connected component is
the image of $\{g\}\times\Omega$ under the uniformization map of Theorem 5.4. Therefore the rule which associates $\chi(\deg(g))^{-1}$ to
the irreducible component corresponding to the double coset represented by
the element $g\in GL_2(\Bbb A_f)$ gives rise to a well-defined $K$-valued
function on the irreducible components of the curve $Y(\goth m)_{F_{\infty}}$.
Actually this function is invariant under the action of the absolute Galois
group of the extension $L$ of $F$ we introduced after Proposition 1.9 hence
the function above is an algebraic cycle on $Y(\goth m)_L$ of co-dimension
zero with coefficients in $K$. For every irreducible component $C$ of $Y(\goth m)_L$ we let $\chi^{-1}(C)$ denote the coefficient of $C$ in this algebraic cycle.
\enddefinition
\definition{Definition 7.6} For every $C$, $D\in\Bbb Z[\Cal V_{\goth m}]_0$ let $\kappa_{\goth m}(C,D)$ denote the element:
$$\epsilon_{\goth m}(C)\otimes\epsilon_{\goth m}(D)\in H^2_{\Cal M}
(Y(\goth m),\Bbb Z(2))$$
where we use the notation of Lemma 5.2. By Proposition 5.5 the pull-back
of $\kappa_{\goth m}(C,D)$ with respect to the uniformization map of Theorem
5.4 is the element of $K_2(GL_2(\Bbb A_f)\times\Omega)$ introduced in Definition 4.9 which is denoted by the same symbol hence our new notation will not cause any confusion. Clearly $\kappa_{\goth m}(C,D)$ is linear in the variables
$C$ and $D$. Let the same symbol denote by abuse of notation the unique
$\Delta$-bilinear extension:
$$\kappa_{\goth m}(\cdot,\cdot):\Delta[\Cal V_{\goth m}]_0\times
\Delta[\Cal V_{\goth m}]_0\rightarrow H^2_{\Cal M}(Y(\goth m),\Delta(2))$$
of this pairing. Let $\kappa_{\goth m}(\chi)$ denote the
unique element of $H^2_{\Cal M}(Y(\goth m)_L,\Delta(2))$ whose restriction
to every irreducible component $C$ of $Y(\goth m)_L$ is $\chi^{-1}(C)
\kappa_{\goth m}(\chi_1,\chi_2)|_C$.
\enddefinition
\proclaim{Theorem 7.7} We have
$$\langle\{\kappa_{\goth m}(\chi)\},\phi_E\rangle
=b(E,\goth m)L(E,q^{-1})\Cal L_{\goth m}(E,\chi)'$$
in $F^*_{\infty}\otimes K$ where $b(E,\goth m)\in K^*$.
\endproclaim
\definition{Proof} The Hecke eigenform $\phi_E$ is locally constant and has compact support as a function on $GL_2(F)\backslash GL_2(\Bbb A)$ hence it takes only finitely many values. In particular there is a positive $n\in\Bbb N$ such that $n\phi_E$ takes integer values. Let
$C$ and $D\in\Delta[\Cal V_{\goth m}]_0$ be two functions such that
the function $\chi(\det(g_f))^{-1}E_{\goth m}^G(C,D,\cdot,x,y)$
is right $Z(\Bbb A)\Bbb K(\goth m\infty)$-invariant. Then the integral:
$$P^G_E(C,D,x,y)=n\langle\chi(\det(g_f))^{-1}E_{\goth m}^G
(C,D,\cdot,x,y),
\phi_E\rangle\in\Delta[[x,y]](x^{-1},y^{-1})$$
is well-defined and it is in fact an element of
$\Delta[x,y,x^{-1},y^{-1}]$ according to Proposition 3.5. Therefore
we may evaluate $P^G_E(C,D,x,y)$ at $x=y=1$. As we already noted at the end of the proof of the limit formula 4.10 we have:
$$E_{\goth m}^G(C,D,g)=(-1)^GE_{\goth m}^G(D,C,g\Pi)$$
for every $g\in GL_2(\Bbb A)$ (using the notation of {\it loc. cit.}). Therefore we get
the equality $P^G_E(C,D,1,1)=-(-1)^GP^G_E(D,C,1,1)$ because
$\phi_E$ is harmonic and the Petersson-product is translation-invariant. The elements $P^G_E(C,D,1,1)$ satisfy the obvious compatibility: let $P_E(C,D)$ denote their limit. Then $P_E(C,D)\in\Delta[[G_{\infty}]]$ lies in $I$ by
Proposition 3.8. Moreover we have:
$$2P_E(\chi_1,\chi_2)'=P_E(\chi_1,\chi_2)'/P_E(\chi_2,\chi_1)'=
\langle\{\kappa_{\goth m}(\chi)\},\phi_E\rangle^n\in G_{\infty}
\otimes K$$
by the Kronecker limit formula 4.10 using the notation we introduced in Definition 6.1. Therefore we get that
$$\langle\{\kappa_{\goth m}(\chi)\},\phi_E\rangle=
{q-1\over2}\mu(\goth m)
L_{\goth m}(E,q^{-1})\Cal L_{\goth m}(E,\chi)'$$
using Theorem 7.4. Since $L_{\goth m}(E,q^{-1})=a(E,\goth m)L(E,q^{-1})$ for some $a(E,\goth m)\in\Bbb Q^*$ the claim follows.\ $\square$
\enddefinition
The function field analogue of the Shimura-Taniyama-Weil conjecture claims the following:
\proclaim{Theorem 7.8} There is a non-trivial map $\pi:X_0(\goth n)\rightarrow
E$ defined over $F$.
\endproclaim
\definition{Proof} Although this theorem is certainly very well known and
have been stated in the literature several times already, in some cases
with an indication of proof, its complete proof have not been written yet,
which we will present now for the sake of record. Let $l$ be a prime
different from $p$ and let $V_{\overline F}$ denote the
base change of any algebraic variety $V$ over $F$ to the separable closure
$\overline F$ of the field $F$. The Gal$(\overline F|F)$-module
$H^1(E_{\overline F}, \Bbb Q_l)$ is absolutely irreducible, because the curve $E$ is not isotrivial. By the global Langlands correspondence for function fields (see [17], proved in this case in [3] already) there is a corresponding
cuspidal automorphic representation $\pi$ of $GL_2(\Bbb A)$. Let
$\omega$ denote the gr\"ossencharacter of $F$ which assigns to each
id\'ele its normalized absolute value. By the compatibility of the local
and global Langlands correspondences the $\infty$-adic component of
$\pi\otimes\omega^{-1}$ is isomorphic to the Steinberg representation.
Also the conductor of $\pi$ is $\goth n\infty$, so there is a non-zero
automorphic form $\phi$ of level $\goth n\infty$ and trivial central
character which is an element of $\pi\otimes\omega^{-1}$. By the above
$\phi$ is also harmonic, so by the main theorem of [4] there is
an absolutely irreducible Gal$(\overline F|F)$ submodule of
$H^1(X_0(\goth n)_{\overline F},\Bbb Q_l)$ corresponding to the
representation $\pi$. This representation must be isomorphic to $H^1
(E_{\overline F},\Bbb Q_l)$, because the Langlands correspondence is a
bijection. By Zarhin's theorem (see [26] and [27]) there is a
homomorphism from the Jacobian of $X_0(\goth n)$ onto $E$ which induces
this isomorphism. We get the map of the claim by composing the map above
with a finite-to-one map from $X_0(\goth n)$ into its Jacobian.\ $\square$
\enddefinition
Our next goal is to give an explicit description of the relation between the modular parameterization of the elliptic curve $E$ in the theorem above and the normalized Hecke eigenform attached to $E$, due to Gekeler and Reversat
[9]. 
\definition{Definition 7.9} Let $\deg(u):GL_2(F_{\infty})
\rightarrow\Bbb Z$ denote the unique function for every holomorphic function
$u:\Omega\rightarrow\Bbb C^*_{\infty}$ such that the regulator
$\{c\otimes f\}$ introduced in Definition 4.5 is equal to
$c^{\deg(u)}$ for every $c\in\Bbb C^*_{\infty}$. Then $\deg(u)$ is just the
van der Put logarithmic derivative of $u$ introduced in [5]. Similarly to
the notation we introduced in Definition 4.7 let $\deg(u):GL_2(\Bbb A)
\rightarrow\Bbb Z$ be the function given by the formula
$\deg(g_f,g_{\infty})=\deg(u(g_f,\cdot))(g_{\infty})$ for each $g\in GL_2(\Bbb A_f)$ if $u:GL_2(\Bbb A_f)\times\Omega\rightarrow\Bbb C_{\infty}^*$ is
ho\-lo\-mor\-phic in the second variable. Recall that $\theta:\Bbb C^*_{\infty}
\rightarrow E(\Bbb C_{\infty})$ denotes the Tate uniformization of $E$.
A theta function attached to $E$ (and the modular parameterization $\pi$)
is a function $u_E:GL_2(\Bbb A_f)\times\Omega\rightarrow\Bbb C_{\infty}^*$
ho\-lo\-mor\-phic in the second variable for each $g\in GL_2(\Bbb A_f)$
if it satisfies the following properties:
\roster
\item"$(a)$" we have $u_E(gk,z)=
u_E(g,z)$ for each $g\in GL_2(\Bbb A_f)$, $z\in\Omega$ and $k\in\Bbb
K_0(\goth n)\cap GL_2(\Bbb A_f)$,
\item"$(b)$" the harmonic cochain $\deg(u_E)$ is $c_E\phi_E$, where
$c_E$ is a positive integer.
\item"$(c)$" the diagram:
$$\CD GL_2(\Bbb A_f)\times\Omega@>>u_E>\Bbb C^*_{\infty}\\
@VVV@VV\theta V\\
Y_0(\goth n)@>>\pi>E(\Bbb C_{\infty})
\endCD$$
is commutative where the vertical map on the left is the uniformization map
mentioned in Notation 5.10.
\endroster
\enddefinition
\proclaim{Theorem (Gekeler-Reversat) 7.10} There is a theta function
attached to $E$.
\endproclaim
\definition{Proof} See [9], Section 9.5, pages 86-88.\ $\square$
\enddefinition
\definition{Proof of Theorem 1.10} Let $\kappa_{\goth m,\goth n}
(\chi)\in H^2_{\Cal M}(Y_0(\goth n)_L,K(2))$ denote the push-forward of the element $\kappa_{\goth m}(\chi)$ with respect to the map $Y(\goth m)
\rightarrow Y_0(\goth m)$ induced the the forgetful map between the
functors represented by these moduli curves. Then we have $\langle\{\kappa_{\goth m}(\chi)\},\psi\rangle=
\langle\{\kappa_{\goth m,\goth n}(\chi)\},\psi\rangle$
for every $\psi\in\Cal H_0(\goth n,\Bbb Q)$ by the invariance theorem (Theorem 3.11 of [22]). Moreover there is a $\kappa_{\goth m,\goth n}(\chi)'\in H^2_{\Cal M}(X_0(\goth n)_L,K(2))$ such that $\langle\{\kappa_{\goth m,\goth n}(\chi)\},
\psi\rangle=\langle\{\kappa_{\goth m,\goth n}(\chi)'\},\psi\rangle$
for every $\psi$ as above by Lemma 5.15. Let $C\subset X_0(\goth n)\times
X_0(\goth n)$ denote the correspondence which is the composition of
the uniformization map $\pi:X_0(\goth n)\rightarrow E$ of Theorem 7.8 and its graph $\Gamma(\pi)\subset E\times X_0(\goth n)$ considered
as a correspondence from $E$ to $X_0(\goth n)$. Then the endomorphism $J(C):J_0(\goth n)\rightarrow J_0(\goth n)$ induced by $C$ is equal to $d(E)T_E$ where $d(E)$ is a non-zero rational number and $T_E\in\Bbb T(\goth n)\otimes\Bbb Q$ is a 
projection operator. Moreover we have $T_E(\phi_E)=\phi_E$ 
therefore
$$\split b(E,\goth m)L(E,q^{-1})\Cal L_{\goth m}(E,\chi)'=&
\langle\{\kappa_{\goth m}(\chi)\},
\phi_E\rangle\\=&\langle\{\kappa_{\goth m,\goth n}
(\chi)'\},\phi_E\rangle\\=&
\langle\{\kappa_{\goth m,\goth n}(\chi)'\},T_E(\phi_E)\rangle\\=&
\langle T_E\{\kappa_{\goth m,\goth n}(\chi)'\},\phi_E\rangle\\=&
d(E)^{-1}\langle\{C_*(\kappa_{\goth m,\goth n}(\chi)')\},\phi_E\rangle\endsplit$$
using the self-adjointness of the operator $T_E$ with respect to
the Petersson product in the fourth equation and Lemma 5.12 in the
last equation. By the invariance theorem (Theorem 3.11 of [22]) the
harmonic form $\{C_*(\kappa_{\goth m,\goth n}(\chi)')\}$ is equal to
$\{\pi_*(\kappa_{\goth m,\goth n}(\chi)')\}\deg(u_E)$ where
$\pi_*(\kappa_{\goth m,\goth n}(\chi)')\in H^2_{\Cal M}(E_L,K(2))$
is the push-forward of $\kappa_{\goth m,\goth n}(\chi)'$
with respect to the unifomization $\pi:X_0(\goth n)\rightarrow E$.
Hence we have
$$L(E,q^{-1})\Cal L_{\goth m}(E,\chi)'=
{c(E)\over b(E,\goth m)d(E)}\{\pi_*(
\kappa_{\goth m,\goth n}(\chi)')\}
\langle\phi_E,\phi_E\rangle$$
by the definition of theta functions. As the Petersson product is
positive definite restricted to $\Cal H_0(\goth n,\Bbb Q)$ the claim is now obvious.\ $\square$
\enddefinition

\heading 8. The action of certain correspondences on $K_2$
\endheading

In this chapter the notation used will be somewhat independent of the one used in the rest of the paper. 
\definition{Notation 8.1} Let $l$ be a prime number and for every scheme $S$ on which $l$ is invertible let
$$c_{2,2}:H^2_{\Cal M}(S,\Bbb Q(2))
\rightarrow H^2_{et}(S,\Bbb Q_l (2))$$
denote the \'etale Chern class map. Let $L$ be a field complete with respect to a discrete valuation and let $\Cal O$ denote its valuation ring. Assume that the residue field of $\Cal O$ is a finite field of characteristic $p\neq l$. Let $\goth X\rightarrow\text{Spec}(\Cal O)$ be a flat, regular, proper and semi-stable scheme over Spec$(\Cal O)$ such that its generic fiber $X$ is a smooth, geometrically irreducible curve over Spec$(L)$. Let $Y$ denote the special fiber of $\goth X$ and let
$$\partial:H^2_{\Cal M}(X,\Bbb Q(2))
\rightarrow H^1_{\Cal M}(Y,\Bbb Q(1))$$ 
denote the boundary map furnished by the localization sequence for the pair $(\goth X,Y)$. 
\enddefinition
\proclaim{Lemma 8.2} For every element $k\in H^2_{\Cal M}(X,\Bbb Q(2))$ such that $c_{2,2}(x) = 0$ we have $\partial(k)=0$.
\endproclaim 
\definition{Proof} Let $R$ denote the ring of global sections of the sheaf of total quotient rings of  $\Cal O_Y$. Since the residue field of any closed point $y$ of $Y$ is a finite field, the homomorphism
$$j^*:H^1_{\Cal M}(Y,\Bbb Q(1))\rightarrow H^1_{\Cal M}(\text{Spec}(R),\Bbb Q(1))=R^*\otimes\Bbb Q$$
induced by the natural map $j:\text{Spec}(R)\rightarrow Y$ is injective. Let
$$c_{1,1}:R^*\otimes\Bbb Q\rightarrow H^1_{et}(\text{Spec}(R),\Bbb Q (1))$$
be the connecting homomorphism of the long cohomological exact sequence attached to Kummer's short exact sequence. Let $\pi:\widetilde Y\rightarrow Y$ be the normalization of $Y$, and let Div$(\widetilde Y)$ denote the group  of divisors on $\widetilde Y$. Since the homomorphism $R^*\rightarrow\text{Div}(\widetilde Y)$ which assigns to every $r\in R^*$ the divisor of $\pi^*(r)$ has a finite kernel and the group Div$(\widetilde Y)$ is a free abelian group, the intersection $\bigcap_{n\in\Bbb N}(R^*)^{l^ n}$ is finite. Hence the homomorphism $c_{1,1}$ is injective. Therefore for every $k\in H^2_{\Cal M}(X,\Bbb Q(2))$ we have $\partial(k)=0$ if the equation $c_{1,1}\circ j_*\circ\partial(k)=0$ holds. Let $K$ be the function field of the curve $X$ and let $i:\text{Spec}(K)\rightarrow X$ be the generic point. The claim now follows form the fact that the diagrams: 
$$\CD H^2_{\Cal M}(X,\Bbb Q(2))@>i^*>>
H^2_{\Cal M}(\text{Spec}(K),\Bbb Q(2))@>\partial>>
R^*\otimes{\Bbb Q}\\
@Vc_{2,2}VV
@Vc_{2,2}VV@Vc_{1,1}VV\\
H^2_{et}(X,\Bbb Q _l(2))@>
i^*>>
H^2_{et}(\text{Spec}(K),\Bbb Q _l(2))@>
\partial>>
H^1_{et}(\text{Spec}(R),\Bbb Q _l(1))\endCD$$
and
$$\CD H^2_{\Cal M}(X,\Bbb Q(2))@>i^*>>
H^2_{\Cal M}(\text{Spec}(K),\Bbb Q(2))\\
@V\partial VV@V\partial VV\\
H^1_{\Cal M}(Y,\Bbb Q(1))@>
j^*>>
R^*\otimes{\Bbb Q}\endCD$$
are commutative, where the symbol $\partial$ denote the respective localisation map everywhere in the diagrams.\ $\square$
\enddefinition
\definition{Notation 8.3} For every smooth, projective, geometrically irreducible curve $Z$ defined over a field $K$ let Jac$(Z)$ denote the Jacobian of $Z$, as usual. Moreover for every correspondence $C\subset Z\times Z$ let
$$C_*:H^2_{\Cal M}(Z,\Bbb Q(2))\rightarrow 
H^2_{\Cal M}(Z,\Bbb Q(2))$$
and
$$J(C):\text{Jac}(Z)\rightarrow\text{Jac}(Z)$$
denote the endomorphisms induced by $C$ on $H^2_{\Cal M}(Z,\Bbb Q(2))$ and Jac$(Z)$, respectively. Let $L$ and $X$ be as in Notation 8.1 and let $C\subset X\times X$ be a correspondence.
\enddefinition
\proclaim{Lemma 8.4} We have $c_{2,2}(C_*(k))=0$ for every
$k\in H^2_{\Cal M}(X,\Bbb Q(2))$ if the endomorphism
$J(C)$ is zero.
\endproclaim 
\definition{Proof} Let $\overline X$ denote the base change of $X$ to the separable closure $\overline L$ of $L$. Note that there is a Hoschschield-Serre spectral sequence:
$$H^i(\text{Gal}(\overline L|L),H^j_{et}(\overline X,\Bbb Q_l(2)))\Rightarrow H^{i+j}_{et}(X,\Bbb Q_l(2)).$$
Because $H^0_{et}(\overline X,\Bbb Q_l(2))=\Bbb Q_l(2)$ and $H^2_{et}(\overline X,\Bbb Q_l(2))=\Bbb Q_l(1)$ we have
$$H^2(\text{Gal}(\overline L|L),H^0_{et}(\overline X,\Bbb Q_l(2)))=0=
H^0(\text{Gal}(\overline L|L),H^2_{et}(\overline X,\Bbb Q_l(2)))$$
by local class field theory. In particular $E^{2,0}_{\infty}=E^{0,2}_{\infty}=0$ for the spectral sequence mentioned above. Moreover $H^3(\text{Gal}(\overline L|L),M)=0$ for every Gal$(\overline L|L)$-module $M$ hence $E^{3,0}_2=0$. Therefore we have $E^{1,1}_{\infty}=E^{1,1}_2=H^1(\text{Gal}(\overline L|L),H^1(\overline X,\Bbb Q_l(2)))$ and so there is an isomorphism:
$$\iota_X:H^2_{et}(X,\Bbb Q_l(2))\rightarrow H^1(\text{Gal}(\overline L|L),H^1(\overline X,\Bbb Q_l(2))).$$
Let 
$$T_l(C):H^1(\overline X,\Bbb Q_l(2)))
\rightarrow H^1(\overline X,\Bbb Q_l(2)))$$
be the endomorphism of $H^1(\overline X,\Bbb Q_l(2)))$ induced by $C$. The map $T_l(C)$ induces a homomorpism on cohomology:
$$T_l(C)_*:H^1(\text{Gal}(\overline L|L),H^1(\overline X,\Bbb Q_l(2)))
\rightarrow
H^1(\text{Gal}(\overline L|L),H^1(\overline X,\Bbb Q_l(2)))$$
by functoriality. Then we have the following commutative diagram:
$$\CD H^2_{\Cal M}(X,\Bbb Q(2))@>c_{2,2}>>
H^2_{et}(X,\Bbb Q_l(2))@>\iota_X>>
H^1(\text{Gal}(\overline L|L),H^1(\overline X,\Bbb Q_l(2)))\\
@VC_*VV
@VC_*VV@VT_l(C)_*VV\\
H^2_{\Cal M}(X,\Bbb Q(2))@>c_{2,2}>>
H^2_{et}(X,\Bbb Q_l(2))@>\iota_X>>
H^1(\text{Gal}(\overline L|L),H^1(\overline X,\Bbb Q_l(2))).\endCD$$
Since we have $H^1(\overline X,\Bbb Q_l(2))=H^1(\text{Jac}(\overline X),\Bbb Q_l(2))$ the map $T_l(C)$ is zero when the endomorphism $J(C)$ is. In this case $T_l(C)_*$ is also zero, so the claim is now clear.\ $\square$
\enddefinition
\definition{Notation 8.5} As in the introduction, let $F$ denote the function field of $X$, where the latter is a geometrically connected smooth projective curve defined over the finite field $\Bbb F_q$ of characteristic $p$. Let $Z$ be a smooth, projective, geometrically irreducible curve defined over $F$ and let $C\subset Z\times Z$ be a correspondence. Assume that $Z$ has a flat, regular, proper and semi-stable model $\goth Z\rightarrow X$ over $X$.
\enddefinition
\proclaim{Proposition 8.6} Assume that the endomorphism $J(C)$ is zero. Then for every $k\in H^2_{\Cal M}(Z,\Bbb Q(2))$ the element $C_*(k)\in H^2_{\Cal M}(Z,\Bbb Q(2))$ lies in the image of the natural map $H^2_{\Cal M}(\goth Z,\Bbb Q(2))\rightarrow H^2_{\Cal M}(Z,\Bbb Q(2))$.
\endproclaim 
\definition{Proof} For every closed point $x$ of $X$ let $\goth Z_x$ denote the fiber of $\goth Z$ at $x$. By the exactness of the localisation sequence it will be sufficient to show that the image of $C_*(k)$ under the boundary map:
$$\partial:H^2_{\Cal M}(Z,\Bbb Q(2))
\rightarrow H^1_{\Cal M}(\goth Z_x,\Bbb Q(1))$$
is zero for every $x$ as above. But this follows at once from Lemmas 8.2 and 8.4 applied to the base change of $Z$ to the completion of $F$ with respect to $x$.\ $\square$
\enddefinition

\enddocument